\documentclass[12pt]{article}
\usepackage{amsfonts, amsmath, amssymb, latexsym, eucal, amscd} 

\usepackage[dvips]{pict2e}

\usepackage{amsthm}
\usepackage{amscd}

\usepackage[dvipdfmx]{graphics}
\usepackage{cite}

\newtheorem{definition}{Definition}[section]
\newtheorem{theorem}[definition]{Theorem}
\newtheorem{lemma}[definition]{Lemma}
\newtheorem{corollary}[definition]{Corollary}
\newtheorem{remark}[definition]{Remark}
\newtheorem{example}[definition]{Example}

\newtheorem{problem}[definition]{Problem}

\newtheorem{assumption}[definition]{Assumption}
\newtheorem{proposition}[definition]{Proposition}

\begin{document} 

\title{\bf 
2-Homogeneous bipartite distance-regular graphs  and the quantum group $U^\prime_q(\mathfrak{so}_6)$
}
\author{
Paul Terwilliger
}
\date{}
\maketitle
\begin{abstract}
We consider a  2-homogeneous bipartite distance-regular graph
$\Gamma$ with diameter $D \geq 3$. We assume that
$\Gamma$ is not a hypercube nor a cycle. 
We fix a $Q$-polynomial ordering of the primitive idempotents of $\Gamma$.
This $Q$-polynomial ordering is described using  a nonzero parameter $q \in \mathbb C$ that is not a root of unity.
We investigate $\Gamma$
using an $S_3$-symmetric approach. In this approach one considers $V^{\otimes 3} = V \otimes V \otimes V$ where $V$
is the standard module of $\Gamma$. We construct a subspace $\Lambda$ of $V^{\otimes 3}$ that has dimension $\binom{D+3}{3}$, together with
six linear maps from $\Lambda$ to $\Lambda$. 
Using these maps we turn $\Lambda$ into an irreducible module for the nonstandard quantum group 
 $U^\prime_q(\mathfrak{so}_6)$ introduced by Gavrilik and Klimyk in 1991.
 \bigskip

 \end{abstract}
 
\section{Introduction} 
The distance-regular graphs have a combinatorial regularity that is often 
studied using algebraic methods \cite{bbit, bannai, bcn,dkt, tSub1, tSub2, tSub3, int}. The hypercubes form an attractive and accessible family of
distance-regular graphs. In \cite{go}, Junie Go turned the standard module $V$ of a $D$-cube into a module for the Lie algebra $\mathfrak{sl}_2(\mathbb C)$.
It is relevant to our story that the Lie algebra  $\mathfrak{sl}_2(\mathbb C)$ is isomorphic to the special orthogonal Lie algebra $\mathfrak{so}_3(\mathbb C)$, see \cite[Section~21.2]{fulton}.
In \cite{martinTer}, Bill Martin and the present author investigated the $D$-cube using the $S_3$-symmetric approach that was introduced in \cite{S3}. In this approach
one considers the vector space  $V^{\otimes 3}=V \otimes V \otimes V$. The authors constructed a subspace $\Lambda$ of $V^{\otimes 3}$ that has dimension
$\binom{D+3}{3}$, along with six linear maps from $\Lambda$ to $\Lambda$. Using these maps the authors turned
 $\Lambda$ into an irreducible module for the Lie algebra $\mathfrak{sl}_4(\mathbb C)$. It is relevant to our story
that the Lie algebra $\mathfrak{sl}_4(\mathbb C)$ is isomorphic to the Lie algebra  $\mathfrak{so}_6(\mathbb C)$, see \cite[Section~21.2]{fulton}.
\medskip

\noindent In \cite{gav2} Gavrilik and Klimyk introduced the nonstandard quantum groups $U'_q(\mathfrak{so}_n)$ $(n\geq 3)$.
The algebra $U'_q(\mathfrak{so}_3)$ or something similar was independently investigated in \cite{fairlie, odesskii, ZheCart}.
The structure and representations of $U'_q(\mathfrak{so}_n)$ were investigated in \cite{gav1, gav2, klimyk, klimykNS, klim1, molev1,noumi2, noumi1, wenzl}.
In \cite{klimyk}, for $q$ not a root of unity the finite-dimensional irreducible $U'_q(\mathfrak{so}_n)$-modules are classified up to isomorphism.
\medskip

\noindent There is a mild generalization of a hypercube called a 2-homogeneous bipartite distance-regular graph. 
In order to clarify what this means, we mention several characterizations of the 2-homogeneous property.
For the rest of this section, let  $\Gamma$ 
denote a bipartite distance-regular graph with diameter $D\geq 3$.  According to \cite[Theorem~42]{curtin},  $\Gamma$ is
2-homogeneous if and only if $\Gamma$ is an antipodal 2-cover and  $Q$-polynomial. Also according to \cite[Theorem~42]{curtin}, $\Gamma$ is 2-homogeneous 
if and only if $\Gamma$ has a $Q$-polynomial ordering of the primitive idempotents that is dual bipartite. In this case,
every $Q$-polynomial ordering of the primitive idempotents  is dual bipartite. For the rest of this section, we assume that $\Gamma$
is  2-homogeneous, but not a hypercube nor a cycle. We fix a $Q$-polynomial ordering of the primitive idempotents of $\Gamma$.
In \cite[Corollaries~36, 43]{curtin} this $Q$-polynomial ordering is described using  a nonzero parameter $q \in \mathbb C$ that is not a root of unity.
In \cite{stronglybalanced} the present author showed that $\Gamma$ has a property called strongly balanced.
In \cite{curtin}, Curtin gave a comprehensive investigation of $\Gamma$ that included many formulas involving $q$. In \cite{curtin2}, Curtin turned the standard module $V$ of $\Gamma$
into a module for $U'_q(\mathfrak{so}_3)$. In \cite{curtin6}, Curtin and Nomura used a weighted adjacency matrix of $\Gamma$
to turn $V$ into a module for the quantum group  $U_q(\mathfrak{sl}_2)$.
\medskip

\noindent In the present paper, we investigate  $\Gamma$ using the $S_3$-symmetric approach.
In rough analogy with \cite{martinTer}, we construct  a subspace $\Lambda$ of $V^{\otimes 3}$ that has dimension
$\binom{D+3}{3}$, along with six linear maps from $\Lambda$ to $\Lambda$. Using these maps we turn  $\Lambda$ into an irreducible module for $U'_q(\mathfrak{so}_6)$. Our main results are
Theorems \ref{thm:main}, \ref{thm:main2}.
\medskip

\noindent The paper is organized as follows. Section 2 contains some preliminaries.
In Section 3 we review some basic concepts and definitions concerning  distance-regular graphs.
In Section 4 we discuss the 2-homogeneous bipartite distance-regular graphs.
In Section 5 we review how these graphs satisfy  the strongly balanced condition.
In Sections 6, 7 we use the strongly balanced condition to establish some combinatorial facts about the graph.
In Section 8 we use these combinatorial facts to construct the vector space $\Lambda$ 
along with six linear maps from  $\Lambda$ to $\Lambda$.
In Section 9 we obtain some relations satisfied by the six maps.
In Section 10 we use the six maps to turn $\Lambda$ into an irreducible module for $U'_q(\mathfrak{so}_6)$.
  Section 11 contains some comments and an open problem. 
Section 12 is an Appendix that contains some technical details.

\section{Preliminaries}

\noindent The following notation and concepts will used throughout the paper.
Recall the natural numbers $\mathbb N = \lbrace 0,1,2,\ldots \rbrace$ and integers $\mathbb Z = \lbrace 0 \pm 1, \pm2, \ldots \rbrace$.
The field of complex numbers is denoted by $\mathbb C$. For $\alpha \in \mathbb C$ let $\overline \alpha$ denote the complex conjugate of $\alpha$.
Every vector space and tensor product that we encounter is understood to be over $\mathbb C$.
Every algebra without the Lie prefix that we encounter, is understood to be associative, over $\mathbb C$, and has a multiplicative
identity. A subalgebra has the same multiplicative identity as the parent algebra. Let $W$ denote a nonzero vector space with finite dimension.
The algebra ${\rm End}(W)$ consists of the $\mathbb C$-linear maps from $W$ to $W$; the algebra product is composition. Let $\mathcal A$ denote an algebra. By an 
{\it automorphism of $\mathcal A$} we mean an algebra isomorphism $\mathcal A \to \mathcal A$.
Let the algebra $\mathcal A^{\rm opp}$ consist of the vector space $\mathcal A$ and the following multiplication.
For $a,b \in \mathcal A$ the product $ab$ (in $\mathcal A^{\rm opp}$) is equal to $ba$ (in $\mathcal A$).
By an {\it antiautomorphism} of $\mathcal A$ we mean an algebra isomorphism $\mathcal A \to \mathcal A^{\rm opp}$.
We will be discussing Hermitian forms. Let us recall the meaning.
\begin{definition}
\label{def:HF} 
 \rm Let $W$ denote a vector space. A {\it Hermitian form on $W$} is a function $\langle \,,\,\rangle: W \times W \to \mathbb C$
such that:
\begin{enumerate}
\item[\rm (i)] $\langle u+v, w \rangle = \langle u, w\rangle + \langle v, w \rangle$ for all $u,v,w \in W$;
\item[\rm (ii)] $\langle \alpha u, v \rangle = \alpha \langle u,v \rangle$ for all $\alpha \in \mathbb C$ and $u,v\in W$;
\item[\rm (iii)] $\langle u,v \rangle = \overline {\langle v,u \rangle}$ for all $u,v \in W$.
\end{enumerate}
\end{definition}
\noindent For a Hermitian form $\langle \,,\,\rangle $ on $W$, we abbreviate $\Vert u \Vert^2 = \langle u,u\rangle$ for all $u \in W$.
\medskip

\section{Distance-regular graphs}
\noindent In this section, we review some definitions and concepts concerning distance-regular graphs. For more information see \cite{bbit, bannai, bcn,dkt, int}.
\medskip

\noindent Let $\Gamma=(X, {\mathcal R})$ denote a finite, undirected, connected graph, without loops or multiple edges, with vertex set $X$, adjacency relation $\mathcal R$,
and path-length distance function $\partial$. To avoid trivialities, we always assume that $\vert X \vert \geq 2$. The positive integer
\begin{align*}
D = {\rm max} \lbrace \partial(x,y) \vert x,y \in X \rbrace
\end{align*}
is called the {\it diameter} of $\Gamma$.
For $0 \leq i \leq D$ and $x \in X$ define the set 
\begin{align*}
\Gamma_i(x) = \lbrace y \in X \vert \partial(x,y)=i \rbrace.
\end{align*}
We abbreviate $\Gamma(x)= \Gamma_1(x)$.
The graph $\Gamma$ is called {\it regular} whenever for $x \in X$ the scalar $k=\vert \Gamma(x) \vert$ is independent of $x$.
In this case, we call $k$ the {\it valency} of $\Gamma$.
The graph $\Gamma$ is called {\it distance-regular} whenever for $0 \leq h,i,j\leq D$ and $x,y \in X$ at distance $\partial(x,y)=h$, the scalar
\begin{align*}
p^h_{i,j} = \vert \Gamma_i(x) \cap \Gamma_j(y) \vert
\end{align*}
is independent of $x,y$ and depends only on $h,i,j$. In this case, the scalars $p^h_{i,j}$ $(0 \leq h,i,j\leq D)$ are called the {\it intersection numbers} of $\Gamma$.
\medskip

\noindent
For the rest of this paper, we assume that $\Gamma$ is distance-regular with diameter $D\geq 3$.
\medskip

\noindent  We comment on the intersection numbers. By construction, $p^h_{i,j} = p^h_{j,i}$ for $0 \leq h,i,j\leq D$. By the triangle inequality, the following hold for $0 \leq h,i,j\leq D$:
\begin{enumerate}
\item[\rm (i)] $p^h_{i,j} =0$ if one of $h,i,j$ is greater than the sum of the other two;
\item[\rm (ii)] $p^h_{i,j} \not=0$ if one of $h,i,j$ is equal to the sum of the other two.
\end{enumerate}
Abbreviate
\begin{align*}
c_i = p^i_{1,i-1} \;(1 \leq i \leq D), \qquad  a_i = p^i_{1,i} \;(0 \leq i \leq D), \qquad b_i = p^i_{1,i+1} \; (0 \leq i \leq D-1)
\end{align*}
and note that $c_1=1$, $a_0=0$. We have $c_i \not=0 $ $(1 \leq i \leq D)$ and $b_i \not=0$ $(0 \leq i \leq D-1)$.
The graph $\Gamma$ is regular with valency $k=b_0$. Moreover
\begin{align*}
c_i + a_i + b_i = k \qquad \qquad (0 \leq i \leq D),
\end{align*}
where $c_0=0$ and $b_D=0$. 
For $0 \leq i \leq D$ define $k_i = p^0_{i,i}$ and note that $k_i = \vert \Gamma_i(x) \vert $ for all $x \in X$.
By construction,
 $\vert X \vert = \sum_{i=0}^D k_i$.
By \cite[p.~247]{bbit} we have
\begin{align*}
k_i = \frac{b_0 b_1 \cdots b_{i-1}}{c_1 c_2 \cdots c_i} \qquad \qquad (0 \leq i \leq D).
\end{align*}
 By \cite[Lemma~3.18]{int} we have
\begin{align*}
       k_h p^h_{i,j} = k_i p^i_{j,h} = k_j p^j_{h,i} \qquad \qquad (0 \leq h,i,j\leq D).
\end{align*}

\noindent Let $V$ denote the vector space with basis $X$. We call $V$ the {\it standard module}. 
 For $x,y \in X$ define $e_{x,y} \in {\rm End}(V)$ that sends $y \mapsto x$ and all other vertices to $0$.
The maps $\lbrace e_{x,y} \rbrace_{x,y \in X}$ form a basis for the vector space ${\rm End}(V)$. 
The transpose map ${\rm t}: {\rm End}(V) \to {\rm End}(V)$ sends $e_{x,y} \mapsto e_{y,x}$ for $x,y \in X$.
The map $\rm t$ is an antiautomorphism of ${\rm End}(V)$.
 We endow $V$ with a Hermitian form $\langle \,,\,\rangle $ with respect to which the basis $X$ is orthonormal.
We have 
\begin{align*}
\langle Bu,v\rangle = \langle u, {\overline B}^{\rm t} v\rangle \qquad \qquad u,v \in V, \qquad B \in {\rm End}(V).
\end{align*}
Define the vector ${\bf 1} \in V$ by
\begin{align}
{\bf 1} =\sum_{x \in X} x.               \label{eq:oneBF}
\end{align}
Define the map $J \in {\rm End}(V)$ by
\begin{align}
J = \sum_{x,y \in X} e_{x,y}.      \label{eq:Jmap}
\end{align}
Note that  $Jx = {\bf 1}$ for all $x \in X$. We have    $\overline J=J=J^{\rm t}$.
\medskip

\noindent We recall the Bose-Mesner algebra of $\Gamma$. 
For $0 \leq i \leq D$ define $A_i \in {\rm End}(V)$ by
\begin{align*}
A_i = \sum_{\stackrel{x,y \in X}{\partial(x,y)=i}}
       e_{x,y}.
\end{align*}
Note that $\overline A_i = A_i=A^{\rm t}_i$. We have
\begin{align*}
A_i x = \sum_{\xi \in \Gamma_i(x)} \xi \qquad \qquad (x \in X).
\end{align*}
Moreover
\begin{align*}
\langle A_i x,y\rangle  = \langle x, A_i y\rangle  =  \begin{cases} 1 & {\mbox{\rm if $\partial(x,y)=i$}}; \\
0, & {\mbox{\rm if $\partial(x,y)\not= i$}}
\end{cases} \qquad \qquad (x,y \in X).
\end{align*}
We call $A_i$ the {\it $i$th distance map} for $\Gamma$.  We abbreviate $A=A_1$ and call $A$ the {\it adjacency map} for $\Gamma$.
We have
\begin{align*}
&A_0 = I, \qquad \qquad J = \sum_{i=0}^D A_i, \\
&A_i A_j = \sum_{h=0}^D p^h_{i,j} A_h \qquad \qquad (0 \leq i,j\leq D).
\end{align*}
The maps $\lbrace A_i \rbrace_{i=0}^D$ form a basis for a commutative subalgebra $M$ of ${\rm End}(V)$. By \cite[Corollary~3.4]{int} the algebra $M$ is generated by $A$.
We call $M$ the {\it Bose-Mesner algebra} of $\Gamma$.
\medskip

\noindent We recall the primitive idempotents of $\Gamma$. 
The map $A$ is real and symmetric, so $A$ is diagonalizable over the real number field. Therefore $M$ has a basis $\lbrace E_i \rbrace_{i=0}^D$ such 
that
\begin{align*}
&E_0 = \vert X \vert^{-1} J, \qquad \qquad \overline E_i = E_i = E_i^{\rm t} \qquad (0 \leq i \leq D), \\
& I = \sum_{i=0}^D E_i, \qquad \qquad  E_i E_j = \delta_{i,j} E_i  \qquad (0 \leq i,j\leq D).
\end{align*}
We call $\lbrace E_i \rbrace_{i=0}^D$ the {\it primitive idempotents} of $M$ (or $\Gamma$).
We have
\begin{align*}
V = \sum_{i=0}^D E_i V \qquad \qquad \hbox{\rm (orthogonal direct sum)}.
\end{align*}
The summands are the eigenspaces of $A$. For $0 \leq i \leq D$ let $\theta_i$ denote the eigenvalue of $A$ associated with $E_iV$.
The scalars $\lbrace \theta_i \rbrace_{i=0}^D$ are real and mutually distinct. Using $AJ=kJ$ we get $\theta_0 = k$.
We have
\begin{align*}
&A = \sum_{i=0}^D \theta_i E_i, \qquad \qquad A E_i = \theta_i E_i = E_i A \qquad  (0 \leq i \leq D), \\
&E_i=\prod_{\stackrel{0 \leq j \leq D}{j \neq i}}
       \frac{A-\theta_jI}{\theta_i-\theta_j}  \qquad \qquad (0 \leq i \leq D).
\end{align*}
\noindent We recall the Krein parameters of $\Gamma$. We turn the vector space ${\rm End}(V)$ into a commutative algebra with the product
\begin{align*}
e_{x,y} \circ e_{x',y'} = \delta_{x,x'} \delta_{y,y'} e_{x,y} \qquad \qquad (x,y, x',y' \in X)
\end{align*}
and multiplicative identity $J$.
We call $\circ$ the {\it Hadamard product}. 
Note that
\begin{align*}
A_i \circ A_j = \delta_{i,j} A_i \qquad \qquad (0 \leq i,j \leq D).
\end{align*}
Thus the Bose-Mesner algebra $M$ is closed under $\circ$. 
Consequently there exist $q^h_{i,j} \in \mathbb C$ $(0 \leq h,i,j\leq D)$  such that
\begin{align*}
E_i \circ E_j = \vert X \vert^{-1} \sum_{h=0}^D q^h_{i,j} E_h \qquad \qquad (0 \leq i,j\leq D).
\end{align*}
By construction, $q^h_{i,j} = q^h_{j,i}$ for $0 \leq h,i,j\leq D$.
By \cite[p.~69]{bannai}, $q^h_{i,j}$ is real and nonnegative for $0 \leq h,i,j\leq D$.
For $0 \leq i \leq D$ define $m_i = q^0_{i,i}$. By  \cite[p.~67]{bannai}, $m_i$ is the dimension of $E_iV$.
By construction, $\vert X \vert=\sum_{i=0}^D m_i$.
 By \cite[p.~67]{bannai} we have
\begin{align*}
       m_h q^h_{i,j} = m_i q^i_{j,h} = m_j q^j_{h,i} \qquad \qquad (0 \leq h,i,j\leq D).
\end{align*}
The scalars $q^h_{i,j}$ $(0 \leq h,i,j \leq D)$ are called the {\it Krein parameters} of $\Gamma$.
\medskip

\noindent We recall the $Q$-polynomial property. The ordering $\lbrace E_i \rbrace_{i=0}^D$ is said to be {\it $Q$-polynomial} whenever
the following hold for $0 \leq h,i,j\leq D$:
\begin{enumerate}
\item[\rm (i)] $q^h_{i,j} =0$ if one of $h,i,j$ is greater than the sum of the other two;
\item[\rm (ii)] $q^h_{i,j} \not=0$ if one of $h,i,j$ is equal to the sum of the other two.
\end{enumerate}
We say that $\Gamma$ is {\it $Q$-polynomial} whenever there exists a $Q$-polynomial ordering of the primitive idempotents of $\Gamma$.
\medskip

\noindent For the rest of this section, we assume that the ordering $\lbrace E_i \rbrace_{i=0}^D$ is $Q$-polynomial.
\medskip

\noindent Abbreviate
\begin{align*}
c^*_i = q^i_{1,i-1} \;(1 \leq i \leq D), \qquad  a^*_i = q^i_{1,i} \;(0 \leq i \leq D), \qquad b^*_i = q^i_{1,i+1} \; (0 \leq i \leq D-1).
\end{align*}
By \cite[p.~67]{bannai} we have $c^*_1=1$, $a^*_0=0$. We have $c^*_i \not=0$ $(1 \leq i \leq D)$ and $b^*_i \not=0$ $(0 \leq i \leq D-1)$.
By \cite[p.~67]{bannai} we have
\begin{align*}
c^*_i + a^*_i + b^*_i = m_1 \qquad \qquad (0 \leq i \leq D),
\end{align*}
where $c^*_0=0$ and $b^*_D=0$. 
By \cite[p.~253]{bbit} we have
\begin{align*}
m_i = \frac{b^*_0 b^*_1 \cdots b^*_{i-1}}{c^*_1 c^*_2 \cdots c^*_i} \qquad \qquad (0 \leq i \leq D).
\end{align*}

\noindent We recall the eigenvalue sequence and dual eigenvalue sequence for the $Q$-polynomial ordering $\lbrace E_i \rbrace_{i=0}^D$.
By \cite[Lemma~19.1]{LSnotes} we have
\begin{align*}
c^*_i \theta_{i-1} + a^*_i \theta_i + b^*_i \theta_{i+1} = \theta^*_1 \theta_i \qquad \qquad (0 \leq i \leq D),
\end{align*}
where $\theta_{-1}$ and $\theta_{D+1}$ denote indeterminates.
The sequence $\lbrace \theta_i \rbrace_{i=0}^D$ is called the {\it eigenvalue sequence} for the 
ordering $\lbrace E_i \rbrace_{i=0}^D$.
For notational convenience, abbreviate $E=E_1$. Since $\lbrace A_i \rbrace_{i=0}^D$ form a basis for $M$, there exist $\theta^*_i \in \mathbb C$ $(0 \leq i \leq D)$
such that
\begin{align*}
E = \vert X \vert^{-1} \sum_{i=0}^D \theta^*_i A_i.
\end{align*}
By \cite[Lemma~11.7]{int} the scalars $\lbrace \theta^*_i \rbrace_{i=0}^D$ are real and mutually distinct. 
By \cite[Lemma~19.1]{LSnotes} we have
\begin{align*}
c_i \theta^*_{i-1} + a_i \theta^*_i + b_i \theta^*_{i+1} = \theta_1 \theta^*_i \qquad \qquad (0 \leq i \leq D),
\end{align*}
where $\theta^*_{-1}$ and $\theta^*_{D+1}$ denote indeterminates.
The sequence $\lbrace \theta^*_i \rbrace_{i=0}^D$ is called the {\it dual eigenvalue sequence} for the 
 ordering $\lbrace E_i \rbrace_{i=0}^D$.

\begin{lemma} 
\label{lem:cosine}
{\rm (See \cite[Proposition~4.4.1]{bcn}.)}
For $x,y \in X$ the following {\rm (i)--(iii)} hold:
\begin{enumerate}
\item[\rm (i)] $
\langle Ex, Ey \rangle = \vert X \vert^{-1} \theta^*_i$ where $i = \partial(x,y)$;
\item[\rm (ii)] $\Vert Ex \Vert^2 = \Vert Ey \Vert^2 = \vert X \vert^{-1} \theta^*_0$;
\item[\rm (iii)]  $\theta^*_i /\theta^*_0$ is the cosine of the angle between $Ex$ and $Ey$.
\end{enumerate}
\end{lemma}

 \begin{corollary} \label{cor:distinct} 
 {\rm (See \cite[Section~4]{NortonBal}.)} 
 For distinct $x,y \in X$ we have $Ex \not=Ey$.
 \end{corollary}
 
 \noindent The graph $\Gamma$ is said to be an {\it antipodal 2-cover}  whenever $k_D=1$.
This occurs if and only if  $b_i = c_{D-i}$  $(0 \leq i \leq D-1)$ if and only if $k_i = k_{D-i} $ $(0 \leq i \leq D)$; see \cite[Lemma~40]{curtin}.
 As we consider additional consequences of Lemma \ref{lem:cosine}, we will treat separately the case in which
 $\Gamma$ is an antipodal 2-cover.
 
 \begin{lemma} \label{lem:NotAntip} 
  {\rm (See \cite[Section~4]{NortonBal}.)} 
 Assume that $\Gamma$ is not an antipodal 2-cover. Then the following hold:
 \begin{enumerate}
 \item[\rm (i)] $\theta^*_0 > \theta^*_i > - \theta^*_0 \quad (1 \leq i \leq D)$;
 \item[\rm (ii)]for distinct $x,y \in X$ the vectors $Ex, Ey$ are linearly independent.
 \end{enumerate}
 \end{lemma}
 
  \begin{lemma} \label{lem:Antip} 
   {\rm (See \cite[Section~4]{NortonBal}.)} 
   Assume that $\Gamma$ is an antipodal 2-cover. Then the following hold:
 \begin{enumerate}
 \item[\rm (i)] $\theta^*_0 > \theta^*_i > - \theta^*_0 \quad (1 \leq i \leq D-1)$ and $\theta^*_D=-\theta^*_0$;
 \item[\rm (ii)]for distinct $x,y \in X$ the vectors $Ex, Ey$ are linearly independent if $\partial(x,y) \not=D$, and  $Ex+Ey=0$ if $\partial(x,y) =D$.
 \end{enumerate}
 \end{lemma}

\begin{lemma} \label{lem:Aix} 
For $x \in X$ and $0 \leq i \leq D$,
\begin{align*}
\sum_{\xi \in \Gamma_i(x)} E \xi = \frac{ k_i \theta^*_i}{\theta^*_0} Ex.
\end{align*}
\end{lemma}
\begin{proof} We have $AE=\theta_1 E$. Since $A$ generates $M$  and $A_i \in M$, there exists a polynomial $f_i$ in one variable such that $A_i = f_i(A)$. Note that
\begin{align*}
\sum_{\xi \in \Gamma_i(x)} E \xi = E A_i x = A_i Ex = f_i(A) Ex = f_i(\theta_1) Ex.
\end{align*}
In this equation we take the inner product of each side with $Ex$ and evaluate the results using Lemma \ref{lem:cosine}; this yields $k_i \theta^*_i = f_i(\theta_1) \theta^*_0$. Therefore $f_i(\theta_1) = k_i \theta^*_i /\theta^*_0$
and the result follows.
\end{proof}

\section{$2$-Homogeneous bipartite distance-regular graphs}

We continue to discuss the distance-regular graph $\Gamma$ with diameter $D\geq 3$.
In this section, we recall the 2-homogeneous bipartite property.
The 2-homogeneous property was introduced by Kazumasa Nomura \cite{2hom}.
In  \cite{curtin, curtin2} Brian Curtin gave a comprehensive treatment of  the case in which $\Gamma$ is 2-homogeneous and bipartite.
\medskip

\noindent The graph $\Gamma$ is called {\it bipartite} whenever $a_i =0$ for $0 \leq i \leq D$.

\begin{definition}
\label{def:2HBip}
\rm
(See \cite[Theorem~42]{curtin}.)
Assume that $\Gamma$ is bipartite. 
Then  $\Gamma$ is said to be {\it $2$-homogeneous} whenever both:
\begin{enumerate}
\item[\rm (i)] $\Gamma$ is an antipodal $2$-cover;
\item[\rm (ii)] $\Gamma$ is $Q$-polynomial.
\end{enumerate}
\end{definition}

\noindent A given $Q$-polynomial ordering $\lbrace E_i \rbrace_{i=0}^D$ is called {\it dual bipartite} whenever
$a^*_i =0$ for $0 \leq i \leq D$.

\begin{lemma}
{\rm (See \cite[Theorem~42]{curtin}.)}
Assume that $\Gamma$ is bipartite.
Then the following are equivalent:
\begin{enumerate}
\item[\rm (i)] $\Gamma$ is $2$-homogeneous;
\item[\rm (ii)] $\Gamma$ has at least one $Q$-polynomial ordering of the primitive idempotents that is dual bipartite.
\end{enumerate}
Assume that {\rm (i), (ii)} hold. Then every  $Q$-polynomial ordering of the primitive idempotents is dual bipartite.
\end{lemma}

\begin{lemma}
{\rm (See \cite[Corollary~43]{curtin} and \cite[Lemma~10.2, Proposition~10.4]{smLP}.)}
 Assume that $\Gamma$ is 2-homogeneous bipartite, and let $\lbrace E_i\rbrace_{i=0}^D$ denote a $Q$-polynomial ordering
of the primitive idempotents of $\Gamma$. Then this ordering is formally self-dual in the sense of \cite[p.~49]{bcn}. In particular:
\begin{align*}
\theta_i &= \theta^*_i \qquad \qquad (0 \leq i \leq D), \\
 p^h_{i,j} &= q^h_{i,j} \qquad \qquad (0 \leq h,i,j\leq D), \\
                  k_i &= m_i \qquad \qquad (0 \leq i \leq D).  
\end{align*}
\end{lemma}

\begin{example}\rm (See \cite[Theorem~1.2]{2hom2}.)
Assume that $\Gamma$ is a hypercube $H(D,2)$ or a $2D$-cycle.
Then $\Gamma$ is 2-homogeneous bipartite.
\end{example}

\begin{lemma} 
\label{lem:q}
{\rm (See \cite[Corollaries~36,~43]{curtin}.)}
Assume that $\Gamma$ is 2-homogeneous bipartite, but not a hypercube nor a cycle.
Let $\lbrace E_i\rbrace_{i=0}^D$ denote a $Q$-polynomial ordering
of the primitive idempotents of $\Gamma$. Then there exists a nonzero $q \in \mathbb C$ that is not a root of unity such that:
\begin{align}
&\theta_i = \theta^*_i = H \sqrt{-1} \bigl(q^{D-2i}- q^{2i-D} \bigr),               \label{eq:thsi} \\
&c_i =c^*_i = H \sqrt{-1} \,\frac{q^{2i}-q^{-2i}}{q^{D-2i}+q^{2i-D}},                     \label{eq:ci}\\
&b_i =b^*_i = H \sqrt{-1}\, \frac{q^{2D-2i}-q^{2i-2D}}{q^{D-2i}+q^{2i-D}}             \label{eq:bi}
\end{align}
for $0 \leq i \leq D$, where 
\begin{align}
H = \sqrt{-1}\, \frac{q^{D-2}+q^{2-D}}{q^{-2}-q^2}.            \label{eq:Hval}
\end{align}
\end{lemma}

\begin{remark} \rm In \cite[Theorem~1.2]{2hom2} Nomura gives a characterization of the 2-homogeneous bipartite distance-regular graphs. This characterization 
 is not used in the present paper.
\end{remark}

\section{The strongly balanced condition}

\noindent From now until the end of Section 9, the following assumption is in effect.

\begin{assumption}
\label{ASSUME} \rm
The graph $\Gamma=(X,\mathcal R)$ is distance-regular with diameter $D\geq 3$.
The graph $\Gamma$ is $2$-homogeneous bipartite, but not a hypercube nor a cycle.
The sequence $\lbrace E_i \rbrace_{i=0}^D$ is a $Q$-polynomial ordering of the primitive idempotents of $\Gamma$.
The corresponding eigenvalue (resp. dual eigenvalue) sequence is  denoted $\lbrace \theta_i \rbrace_{i=0}^D$ (resp. $\lbrace \theta^*_i \rbrace_{i=0}^D$).
The scalar $q$ is from Lemma \ref{lem:q}.
\end{assumption}

\noindent  In this section we review the strongly balanced condition. One version of
this condition is given in the next result.

\begin{proposition} \label{prop:SBC} {\rm (See \cite[Theorems~1, 3]{stronglybalanced}.)} For $x,y \in X$ and $0 \leq i,j\leq D$,
\begin{align*}
\sum_{\xi \in \Gamma_i(x) \cap \Gamma_j(y)} E \xi \in {\rm Span} \lbrace Ex, Ey\rbrace.
\end{align*}
\end{proposition}

\noindent Next, we give a more detailed version of Proposition \ref{prop:SBC} that shows the coefficients.
\begin{proposition} \label{cor:SBC} 
Let $1 \leq h \leq D-1$ and $x,y \in X$ at distance $\partial(x,y) = h$. Then for $0 \leq i,j\leq D$,
\begin{align*}
\sum_{\xi \in \Gamma_i(x) \cap \Gamma_j(y)} E \xi = 
p^h_{i,j} \frac{\theta^*_0 \theta^*_i - \theta^*_h \theta^*_j}{ \theta^{* 2}_0 - \theta^{* 2}_h} Ex+ 
p^h_{i,j} \frac{\theta^*_0 \theta^*_j - \theta^*_h \theta^*_i}{ \theta^{* 2}_0 - \theta^{* 2}_h} Ey.
\end{align*}
\end{proposition} 
\begin{proof}  By Proposition \ref{prop:SBC} there exist complex scalars $r^h_{i,j}$ and $s^h_{i,j}$ such that
\begin{align}
\sum_{\xi \in \Gamma_i(x) \cap \Gamma_j(y)} E \xi = r^h_{i,j} Ex + s^h_{i,j} Ey.        \label{eq:RS}
\end{align}
In \eqref{eq:RS} we take the inner product of each side with $Ex$ and evaluate the results using Lemma \ref{lem:cosine}; this yields
\begin{align}
p^h_{i,j} \theta^*_i = r^h_{i,j} \theta^*_0 + s^h_{i,j} \theta^*_h. \label{eq:one}
\end{align}
In \eqref{eq:RS} we take the inner product of each side with $Ey$; this similarly yields
\begin{align}
p^h_{i,j} \theta^*_j = r^h_{i,j} \theta^*_h + s^h_{i,j} \theta^*_0. \label{eq:two}
\end{align}
Using Lemma \ref{lem:Antip}  and linear algebra, we solve \eqref{eq:one}, \eqref{eq:two} to obtain
\begin{align} \label{eq:rsAns}
r^h_{i,j} = p^h_{i,j} \frac{\theta^*_0 \theta^*_i - \theta^*_h \theta^*_j}{ \theta^{* 2}_0 - \theta^{* 2}_h}, \qquad \qquad 
s^h_{i,j} = p^h_{i,j} \frac{\theta^*_0 \theta^*_j - \theta^*_h \theta^*_i}{ \theta^{* 2}_0 - \theta^{* 2}_h}.
\end{align}
Combining \eqref{eq:RS}, \eqref{eq:rsAns} we get the result.
\end{proof} 

\noindent We mention two special cases of Proposition \ref{cor:SBC}.
\begin{corollary}
 \label{cor:SBsc1} 
 Let $1 \leq h \leq D-1$ and $x,y \in X$ at distance $\partial(x,y)=h$. Then
\begin{align*}
\sum_{\xi \in \Gamma(x) \cap \Gamma_{h-1}(y)} E \xi = 
c_h \frac{\theta^*_0 \theta^*_1 - \theta^*_{h-1} \theta^*_h}{ \theta^{* 2}_0 - \theta^{* 2}_h} Ex+ 
c_h \frac{\theta^*_0 \theta^*_{h-1} - \theta^*_1 \theta^*_h}{ \theta^{* 2}_0 - \theta^{* 2}_h} Ey.
\end{align*}
\end{corollary}
\begin{proof} This is Proposition \ref{cor:SBC} with $i=1$ and $j=h-1$.
\end{proof}
\begin{corollary}
\label{cor:SBsc2}
 Let $1 \leq h \leq D-1$ and $x,y \in X$ at distance $\partial(x,y)=h$. Then
\begin{align*}
\sum_{\xi \in \Gamma(x) \cap \Gamma_{h+1}(y)} E \xi = 
b_h \frac{\theta^*_0 \theta^*_1 - \theta^*_{h+1} \theta^*_h}{ \theta^{* 2}_0 - \theta^{* 2}_h} Ex+ 
b_h \frac{\theta^*_0 \theta^*_{h+1} - \theta^*_1 \theta^*_h}{ \theta^{* 2}_0 - \theta^{* 2}_h} Ey.
\end{align*}
\end{corollary}
\begin{proof} This is Proposition \ref{cor:SBC} with $i=1$ and $j=h+1$.
\end{proof}

\noindent For notational convenience, we define $\Gamma_{-1}(x) = \emptyset$ and $\Gamma_{D+1} (x) = \emptyset$ for $x \in X$.

\section{Some combinatorial regularity}

\noindent We continue to discuss the graph $\Gamma$ from Assumption \ref{ASSUME}.
Throughout this section, we fix $x, y, z \in X$ and write
\begin{align*}
h= \partial(y,z), \qquad \qquad i = \partial(z,x), \qquad \qquad j = \partial(x,y).
\end{align*}
\noindent  Since $\Gamma$ is bipartite, for $\xi \in \Gamma(x)$ we have
$\partial(\xi, z) \in \lbrace i-1, i+1\rbrace$ and $\partial(\xi, y) \in \lbrace j-1, j+1\rbrace$. Thus
 the set $\Gamma(x)$ is a disjoint union
of the following four sets (some  might be empty):
\begin{align} \label{eq:sets12}
&\Gamma(x) \cap \Gamma_{i+1} (z) \cap \Gamma_{j+1}(y), \qquad \qquad 
\Gamma(x) \cap \Gamma_{i-1} (z) \cap \Gamma_{j-1}(y),  \\
&\Gamma(x) \cap \Gamma_{i+1} (z) \cap \Gamma_{j-1}(y), \qquad \qquad
\Gamma(x) \cap \Gamma_{i-1} (z) \cap \Gamma_{j+1}(y).
\label{eq:sets34}
\end{align}
\noindent  Our next goal is to compute the cardinalities of the above four sets.
As we will see, these cardinalities  depend only on $h,i,j$ and not on the choice of $x,y,z$.
\medskip

\begin{lemma} \label{lem:easy} We have
\begin{align*}
c_i &=\vert \Gamma(x) \cap \Gamma_{i-1} (z) \cap \Gamma_{j-1}(y) \vert + \vert \Gamma(x) \cap \Gamma_{i-1} (z) \cap \Gamma_{j+1}(y) \vert, \\
c_j &= \vert \Gamma(x) \cap \Gamma_{i-1} (z) \cap \Gamma_{j-1}(y) \vert + \vert \Gamma(x) \cap \Gamma_{i+1} (z) \cap \Gamma_{j-1}(y) \vert,     \\
b_i &= \vert \Gamma(x) \cap \Gamma_{i+1} (z) \cap \Gamma_{j-1}(y) \vert + \vert \Gamma(x) \cap \Gamma_{i+1} (z) \cap \Gamma_{j+1}(y) \vert,       \\
b_j &= \vert \Gamma(x) \cap \Gamma_{i-1} (z) \cap \Gamma_{j+1}(y) \vert + \vert \Gamma(x) \cap \Gamma_{i+1} (z) \cap \Gamma_{j+1}(y) \vert . 
\end{align*}
\end{lemma}
\begin{proof} To obtain the first equation in the lemma statement, observe that $c_i=\vert \Gamma(x) \cap \Gamma_{i-1}(z)\vert $
and $\Gamma(x) \cap \Gamma_{i-1}(z)$ is the disjoint union of the sets on the right in \eqref{eq:sets12}, \eqref{eq:sets34}. The remaining equations in the lemma statement are similarly obtained.
\end{proof}

\begin{lemma} \label{lem:hard} Assume that $1 \leq i,j\leq D-1$. Then
\begin{align*}
&\theta^*_{j-1} \vert \Gamma(x) \cap \Gamma_{i-1} (z) \cap \Gamma_{j-1}(y) \vert  +
\theta^*_{j+1}  \vert \Gamma(x) \cap \Gamma_{i-1} (z) \cap \Gamma_{j+1}(y) \vert  \\
& = 
c_i \theta^*_j \frac{\theta^*_0 \theta^*_1 - \theta^*_{i-1} \theta^*_i}{ \theta^{* 2}_0 - \theta^{* 2}_i} + 
c_i \theta^*_h \frac{\theta^*_0 \theta^*_{i-1} - \theta^*_1 \theta^*_i}{ \theta^{* 2}_0 - \theta^{* 2}_i},
\end{align*}
\begin{align*}
&\theta^*_{i-1} \vert \Gamma(x) \cap \Gamma_{i-1} (z) \cap \Gamma_{j-1}(y) \vert  +
\theta^*_{i+1}  \vert \Gamma(x) \cap \Gamma_{i+1} (z) \cap \Gamma_{j-1}(y) \vert  \\
& = 
c_j \theta^*_i \frac{\theta^*_0 \theta^*_1 - \theta^*_{j-1} \theta^*_j}{ \theta^{* 2}_0 - \theta^{* 2}_j} + 
c_j \theta^*_h \frac{\theta^*_0 \theta^*_{j-1} - \theta^*_1 \theta^*_j}{ \theta^{* 2}_0 - \theta^{* 2}_j}.
\end{align*}
\begin{align*}
&\theta^*_{j-1} \vert \Gamma(x) \cap \Gamma_{i+1} (z) \cap \Gamma_{j-1}(y) \vert  +
\theta^*_{j+1}  \vert \Gamma(x) \cap \Gamma_{i+1} (z) \cap \Gamma_{j+1}(y) \vert  \\
& = 
b_i \theta^*_j \frac{\theta^*_0 \theta^*_1 - \theta^*_{i+1} \theta^*_i}{ \theta^{* 2}_0 - \theta^{* 2}_i} + 
b_i \theta^*_h \frac{\theta^*_0 \theta^*_{i+1} - \theta^*_1 \theta^*_i}{ \theta^{* 2}_0 - \theta^{* 2}_i},
\end{align*}
\begin{align*}
&\theta^*_{i-1} \vert \Gamma(x) \cap \Gamma_{i-1} (z) \cap \Gamma_{j+1}(y) \vert  +
\theta^*_{i+1}  \vert \Gamma(x) \cap \Gamma_{i+1} (z) \cap \Gamma_{j+1}(y) \vert  \\
& = 
b_j \theta^*_i \frac{\theta^*_0 \theta^*_1 - \theta^*_{j+1} \theta^*_j}{ \theta^{* 2}_0 - \theta^{* 2}_j} + 
b_j \theta^*_h \frac{\theta^*_0 \theta^*_{j+1} - \theta^*_1 \theta^*_j}{ \theta^{* 2}_0 - \theta^{* 2}_j}.
\end{align*}
\end{lemma}
\begin{proof} We obtain the first equation in the lemma statement. By Lemma  \ref{cor:SBsc1},
\begin{align}
\label{eq:3step4}
\sum_{\xi \in \Gamma(x) \cap \Gamma_{i-1}(z)} E \xi
& = 
c_i  \frac{\theta^*_0 \theta^*_1 - \theta^*_{i-1} \theta^*_i}{ \theta^{* 2}_0 - \theta^{* 2}_i} Ex+ 
c_i  \frac{\theta^*_0 \theta^*_{i-1} - \theta^*_1 \theta^*_i}{ \theta^{* 2}_0 - \theta^{* 2}_i} Ez.
\end{align}
In the equation \eqref{eq:3step4}, take the inner product of each side with $Ey$ and evaluate the result using
 Lemma \ref{lem:cosine}. This yields the first equation in the lemma statement. The remaining equations
 of the lemma statement are similarly obtained.
\end{proof}

\begin{proposition} \label{thm:2hom} Pick $x, y, z \in X$ and write
\begin{align*}
h= \partial(y,z), \qquad \qquad i = \partial(z,x), \qquad \qquad j = \partial(x,y).
\end{align*}
\noindent Then
\begin{align*}
\vert \Gamma(x) \cap \Gamma_{i+1} (z) \cap \Gamma_{j+1}(y) \vert &= H\sqrt{-1}\,\frac{q^{2D-h-i-j}-q^{h+i+j-2D}}{q^{D-2i}+q^{2i-D}} \, \frac{q^{D+h-i-j}+ q^{i+j-h-D}}{q^{D-2j}+q^{2j-D}}, \\
\vert \Gamma(x) \cap \Gamma_{i-1} (z) \cap \Gamma_{j-1}(y) \vert &= H\sqrt{-1}\,\frac{q^{i+j-h}-q^{h-i-j}}{q^{D-2i}+q^{2i-D}} \, \frac{q^{D-h-i-j}+ q^{h+i+j-D}}{q^{D-2j}+q^{2j-D}}, \\
\vert \Gamma(x) \cap \Gamma_{i+1} (z) \cap \Gamma_{j-1}(y) \vert &= H\sqrt{-1}\,\frac{q^{h+j-i}-q^{i-h-j}}{q^{D-2i}+q^{2i-D}} \, \frac{q^{D+j-h-i}+ q^{h+i-j-D}}{q^{D-2j}+q^{2j-D}}, \\
\vert \Gamma(x) \cap \Gamma_{i-1} (z) \cap \Gamma_{j+1}(y) \vert &= H\sqrt{-1}\,\frac{q^{h+i-j}-q^{j-h-i}}{q^{D-2i}+q^{2i-D}} \, \frac{q^{D+i-h-j}+ q^{h+j-i-D}}{q^{D-2j}+q^{2j-D}}.
\end{align*}
\end{proposition}
\begin{proof} First assume that $1 \leq i,j\leq D-1$. To obtain the first and third equations in the proposition statement, combine the third equation in Lemma \ref{lem:easy}
and the third equation in Lemma \ref{lem:hard}. Evaluate the result using \eqref{eq:thsi} and \eqref{eq:bi}.  To obtain the second and fourth equations in the proposition statement,
combine the first equation in Lemma \ref{lem:easy}  and the first equation in Lemma \ref{lem:hard}.
 Evaluate the result using \eqref{eq:thsi} and \eqref{eq:ci}. 
Next assume that $i=0$ or $i=D$ or $j=0$ or $j=D$. In each case, the four equations in the proposition statement are routinely
checked using \eqref{eq:ci}, \eqref{eq:bi} and the fact that $\Gamma$ is a bipartite antipodal 2-cover.
\end{proof}

%

\section{A change of variables}
\noindent We continue to discuss the graph $\Gamma$ from Assumption \ref{ASSUME}.
 In this section
we describe a change of variables that improves the formulas in Proposition \ref{thm:2hom}.

\begin{definition} \label{def:pf} \rm Let the set $\mathcal P_D$ consist of the 4-tuples of natural numbers $(r,s,t,u)$ such that
$r+s+t+u=D$. An element of $\mathcal P_D$ is called a {\it profile of degree $D$}.
\end{definition}
\noindent Note that 
\begin{align}
\vert \mathcal P_D \vert = \binom{D+3}{3}.                      \label{eq:pfSIZE}
\end{align}

\begin{definition} \label{def:prP} \rm Let the set $\mathcal P'_D$ consist of the $3$-tuples of integers $(h,i,j)$ such that
\begin{align*}
&0 \leq h,i,j\leq D, \qquad  h+i+j \;\hbox{\rm is even}, \qquad  h+i+j \leq 2D, \\
&h \leq i+j, \quad \qquad \qquad i \leq j+h, \quad \qquad \qquad j \leq h+i.
\end{align*}
\end{definition}

\begin{lemma} 
\label{lem:nonzINT}
 {\rm (See \cite[Corollary~4.3.11]{nicholson}.)} For $0 \leq h,i,j\leq D$ the intersection number $p^h_{i,j}$ is nonzero if and only if
$(h,i,j) \in \mathcal P'_D$.
\end{lemma}

\begin{lemma} 
\label{lem:prBIJ} 
There exists a bijection $\mathcal P_D \to \mathcal P'_D$ that sends
\begin{align*}
(r,s,t,u) \mapsto (t+u,u+s,s+t).
\end{align*}
The inverse bijection  $\mathcal P'_D \to \mathcal P_D$ sends
\begin{align*}
(h,i,j) \mapsto \biggl( \frac{2D-h-i-j}{2}, \frac{i+j-h}{2}, \frac{j+h-i}{2}, \frac{h+i-j}{2}\biggr).
\end{align*}
\end{lemma}
\begin{proof} This is routinely  checked.
\end{proof}

\begin{lemma} 
\label{lem:interp}
 For a 3-tuple of vertices $x,y,z$  there exists a unique profile $(r,s,t,u) \in {\mathcal P}_D$ such that
\begin{align*}
\partial(x,y)= s+t, \qquad \quad \partial(y,z)=t+u, \qquad \quad \partial(z,x) = u+s.
\end{align*}
Moreover
\begin{align*}
&r = \frac{ 2D - \partial(x,y) - \partial(y,z)-\partial(z,x)}{2}, \qquad \quad s = \frac{\partial(z,x)+\partial(x,y)-\partial(y,z)}{2}, \\
&t = \frac{\partial(x,y)+\partial(y,z)-\partial(z,x)}{2}, \qquad \quad u = \frac{\partial(y,z)+\partial(z,x)-\partial(x,y)}{2}.
\end{align*}
\end{lemma}
\begin{proof} By Lemmas \ref{lem:nonzINT}, \ref{lem:prBIJ}.
\end{proof}

\begin{definition} 
\label{def:profile} 
\rm Referring to  Lemma \ref{lem:interp}, we call $(r,s,t,u)$  the {\it profile} of $x,y,z$.
\end{definition}

\begin{lemma} 
\label{lem:AView}
For  $0 \leq h,i,j\leq D$ the following {\rm (i)--(iii)}
are equivalent:
\begin{enumerate}
\item[\rm (i)] there exist $x,y,z \in X$  such that
\begin{align*}
h=\partial(y,z), \qquad \qquad i=\partial(z,x), \qquad \qquad
j = \partial(x,y);
\end{align*}
\item[\rm (ii)]  $(h,i,j) \in {\mathcal P}'_D$;
\item[\rm (iii)]  there exists $(r,s,t,u) \in \mathcal P_D$ such that
\begin{align*}
 h=t+u, \qquad \qquad i = u+s, \qquad \qquad j= s+t.
 \end{align*}
\end{enumerate}
\noindent Assume that {\rm (i)--(iii)} hold. Then $(r,s,t,u)$ is the profile of $x,y,z$.
\end{lemma}
\begin{proof} ${\rm (i)} \Leftrightarrow {\rm (ii)}$ By  Definition \ref{def:prP} and Lemma \ref{lem:nonzINT}. \\
${\rm (ii)} \Leftrightarrow {\rm (iii)}$ By Lemma \ref{lem:prBIJ}. \\
The last assertion follows from Lemma \ref{lem:interp} and Definition \ref{def:profile}.    
\end{proof}

\begin{definition} \label{def:C} \rm For $(r,s,t,u ) \in {\mathcal P}_D$  define
\begin{align}
C(r,s,t,u)= H\sqrt{-1}\, \frac{q^{2r}-q^{-2r}}{q^{r+u-s-t}+ q^{s+t-r-u}}\, \frac{q^{t-r-s-u}+q^{r+s+u-t}}{q^{r+s-u-t}+q^{u+t-r-s}} \label{eq:C}
\end{align}
where $H$ is from  \eqref{eq:Hval}.
\end{definition}

\noindent In the next result, we express Proposition  \ref{thm:2hom}  in terms of profiles.

\begin{proposition}
\label{thm:2homALT}
Pick $x, y, z \in X$ and write
\begin{align*}
h= \partial(y,z), \qquad \qquad i = \partial(z,x), \qquad \qquad j = \partial(x,y).
\end{align*}
\noindent Then
\begin{align*}
\vert \Gamma(x) \cap \Gamma_{i+1} (z) \cap \Gamma_{j+1}(y) \vert &= C(r,t,s,u), \\
\vert \Gamma(x) \cap \Gamma_{i-1} (z) \cap \Gamma_{j-1}(y) \vert &= C(s,u,r,t), \\
\vert \Gamma(x) \cap \Gamma_{i+1} (z) \cap \Gamma_{j-1}(y) \vert &= C(t,s,u,r), \\
\vert \Gamma(x) \cap \Gamma_{i-1} (z) \cap \Gamma_{j+1}(y) \vert &= C(u,r,t,s).
\end{align*}
where $(r,s,t,u)$ is the profile of $x,y,z$.
\end{proposition}
\begin{proof} Evaluate Proposition \ref{thm:2hom} using
Lemma  \ref{lem:interp} and Definition \ref{def:C}.
\end{proof}

\section{The graph $\Gamma$ from an $S_3$-symmetric point of view}

\noindent We continue to discuss the graph $\Gamma$ from Assumption \ref{ASSUME}.
 Let $S_3$ denote
the symmetric group on the set $\lbrace 1,2,3\rbrace$.
In this section, we investigate $\Gamma$ from the $S_3$-symmetric point of view introduced in \cite{S3}.
\medskip

\noindent Recall the standard module $V$ of $\Gamma$ from Section 3.

\begin{definition}\label{def:VVV} \rm
Define the vector space $V^{\otimes 3} = V \otimes V \otimes V$ and the set
\begin{align*}
X^{\otimes 3} = \lbrace x \otimes y \otimes z \vert x,y,z \in X\rbrace.
\end{align*}
\end{definition}
\noindent Note that $X^{\otimes 3}$ is a basis for $V^{\otimes 3}$.

\begin{lemma} 
\label{lem:HermV3} 
There exists a unique Hermitian form $\langle \,,\,\rangle$ on $V^{\otimes 3}$ with respect to which $X^{\otimes 3}$ is
orthonormal. For $u\otimes v \otimes w \in V^{\otimes 3}$ and
 $u'\otimes v' \otimes w' \in V^{\otimes 3}$,
 \begin{align*}
\bigl \langle u \otimes v \otimes w, u' \otimes v' \otimes w' \bigr \rangle = \langle u, u' \rangle \langle v, v' \rangle \langle w, w' \rangle.
 \end{align*}
 \end{lemma}
\begin{proof} By linear algebra.
\end{proof}

\noindent Our next goal is to introduce some maps in ${\rm End}(V^{\otimes 3})$, denoted
\begin{align}
A^{(1)}, \quad A^{(2)}, \quad A^{(3)}, \quad A^{*(1)}, \quad A^{*(2)}, \quad A^{*(3)}.
\label{eq:gens}
\end{align}

\begin{definition}
\label{def:maps}
\rm 
We define $A^{(1)}, A^{(2)}, A^{(3)} \in {\rm End}(V^{\otimes 3}) $ as follows.
For $x\otimes y\otimes z \in X^{\otimes 3}$,
\begin{align*}
A^{(1)} (x \otimes y \otimes z) &= \sum_{\xi \in \Gamma(x)} \xi \otimes y \otimes z, \\
A^{(2)} (x \otimes y \otimes z) &= \sum_{\xi \in \Gamma(y)} x \otimes \xi \otimes z, \\
A^{(3)} (x \otimes y \otimes z) &= \sum_{\xi \in \Gamma(z)} x \otimes y \otimes \xi.
\end{align*}
\end{definition}

\noindent The next result is meant to clarify Definition \ref{def:maps}.

\begin{lemma} 
\label{lem:clarify}
For $u \otimes v \otimes w \in V^{\otimes 3}$ we have
\begin{align*}
A^{(1)} (u \otimes v \otimes w) &= (Au) \otimes v \otimes w, \\
A^{(2)} (u \otimes v \otimes w) &= u \otimes (Av) \otimes w, \\
A^{(3)} (u \otimes v \otimes w) &= u \otimes v \otimes (Aw).
\end{align*}
\end{lemma} 
\begin{proof} By the definition of the adjacency map $A$.
\end{proof}

\begin{definition}
\label{def:dualmaps}
\rm 
We define $A^{*(1)}, A^{*(2)}, A^{*(3)} \in {\rm End}(V^{\otimes 3}) $ as follows.
For $x\otimes y\otimes z \in X^{\otimes 3}$,
\begin{align*}
A^{*(1)} (x \otimes y \otimes z) &= \theta^*_{\partial(y,z)} x \otimes y \otimes z, \\
A^{*(2)} (x \otimes y \otimes z) &= \theta^*_{\partial(z,x)} x \otimes y \otimes z, \\
A^{*(3)} (x \otimes y \otimes z) &= \theta^*_{\partial(x,y)} x \otimes y \otimes z.
\end{align*}
\end{definition}

\begin{lemma}
\label{lem:sym}
For $i \in \lbrace 1,2,3\rbrace$ and $u,v \in V^{\otimes 3}$ we have
\begin{align}
\bigl \langle A^{(i)} u,v \bigr \rangle = \bigl \langle u, A^{(i)} v \bigr \rangle,  \qquad \qquad         \bigl \langle A^{*(i)} u,v \bigr \rangle = \bigl \langle u, A^{*(i)} v \bigr \rangle.
\end{align}
\end{lemma}
\begin{proof} By $S_3$-symmetry we may assume without loss of generality that $i=1$.
Also without loss of generality, we may assume that $u,v$ are contained in the basis $X^{\otimes 3}$. Write
$u = x \otimes y \otimes z$ and $v = x' \otimes y' \otimes z'$. By Definition \ref{def:maps},
\begin{align*}
\bigl \langle A^{(1)} (x\otimes y \otimes z), x' \otimes y' \otimes z' \bigr \rangle &=
\begin{cases} 1   & \hbox{\rm if $\partial(x,x')=1$ and $y=y'$ and $z = z'$};\\
                       0   & \hbox{\rm if $\partial(x,x')\not=1$ or $y\not=y' $ or $z \not= z'$}
\end{cases}
\\
&= \bigl \langle  x\otimes y \otimes z, A^{(1)} (x' \otimes y' \otimes z') \bigr \rangle.
\end{align*}
By Definition \ref{def:dualmaps} and since the dual eigenvalues are real,
\begin{align*}
\bigl \langle A^{*(1)} (x\otimes y \otimes z), x' \otimes y' \otimes z' \bigr \rangle &=
\begin{cases}  \theta^*_{\partial(y,z)}
 & \hbox{\rm if $x=x'$ and $y=y'$ and $z = z'$};\\
                    0   & \hbox{\rm if $x\not=x'$ or $y\not=y' $ or $z \not= z'$}
\end{cases}
\\
&= \bigl \langle  x\otimes y \otimes z, A^{*(1)} (x' \otimes y' \otimes z')\bigr \rangle.
\end{align*}
\noindent The result follows.
\end{proof}

\begin{definition}
\label{def:norm} \rm
For a profile $(r,s,t,u) \in {\mathcal P}_D$ let $n(r,s,t,u)$ denote the number of 3-tuples of vertices $x,y,z$ that have profile $(r,s,t,u)$.
\end{definition}

\begin{lemma}
\label{lem:norm}
Pick a profile $(r,s,t,u) \in {\mathcal P}_D$ and write
\begin{align*}
h=t+u, \qquad \quad i = u+s, \qquad \quad j = s+t.
\end{align*}
 Then $n(r,s,t,u)$ is equal to the number of $3$-tuples of vertices $x,y,z$ such that 
\begin{align*}
h= \partial(y,z), \qquad \quad i = \partial(z,x), \qquad \quad j = \partial(x,y).
\end{align*}
Moreover,
\begin{align*}
n(r,s,t,u) = \vert X \vert k_h p^h_{i,j} = \vert X \vert k_i p^i_{j,h} = \vert X \vert k_j p^j_{h,i}\not=0.
\end{align*}
\end{lemma}
\begin{proof} By Lemmas \ref{lem:nonzINT}, \ref{lem:prBIJ}, \ref{lem:interp}   and combinatorial counting.
\end{proof}

\begin{definition}
\label{def:Bcheck}
\rm
For a profile $(r,s,t,u) \in {\mathcal P}_D$ define a vector
\begin{align*}
B (r,s,t,u) = \sum x \otimes y \otimes z,
\end{align*}
where the sum is over the $3$-tuples of vertices $x,y,z$ that have profile $(r,s,t,u)$.
\end{definition}

\begin{example} \label{ex:BDzzz} \rm We have
\begin{align*}
B(D,0,0,0) = \sum_{x \in X}   x \otimes x \otimes x.
\end{align*}
\end{example}

\noindent We remark about the notation.
\begin{remark}
\label{rem:aside}
\rm
(See \cite[Definition~9.9]{S3}.)
For $0 \leq h,i,j\leq D$ define a vector 
\begin{align*}
P_{h,i,j} = \sum x \otimes y \otimes z,
\end{align*}
where the sum is over the $3$-tuples of vertices $x,y,z$ such that
\begin{align*}
h = \partial(y,z), \qquad \quad i = \partial(z,x), \qquad \quad j = \partial(x,y).
\end{align*}
By Lemma \ref{lem:AView} we have the following.
The vector $P_{h,i,j} \not=0$ if and only if $(h,i,j) \in {\mathcal P}'_D$.
In this case,
$P_{h,i,j} = B (r,s,t,u)$ where $(r,s,t,u) \in {\mathcal P}_D$ satisfies
\begin{align*}
h=t+u, \qquad \quad i = u+s, \qquad \quad j = s+t.
\end{align*}
\end{remark}

\begin{lemma}
\label{lem:Borthog} 
The vectors
\begin{align*}
B(r,s,t,u) \qquad \qquad (r,s,t,u) \in {\mathcal P}_D            
\end{align*}
are mutually orthogonal. Moreover
\begin{align*}
\Vert B (r,s,t,u) \Vert^2 = n(r,s,t,u) \qquad \qquad (r,s,t,u) \in {\mathcal P}_D.
\end{align*}
\end{lemma}
\begin{proof} By Definition \ref{def:Bcheck} and since the basis $X^{\otimes 3}$ is orthonormal with respect to $\langle \,,\,\rangle$.
\end{proof}

\begin{definition}
\label{def:Lambda}
\rm Let $\Lambda$ denote the subspace of $V^{\otimes 3}$ with the basis 
\begin{align*}
B(r,s,t,u) \qquad \qquad (r,s,t,u) \in {\mathcal P}_D.           
\end{align*}
\end{definition}

\begin{lemma} 
\label{lem:LambdaDim} The vector space $\Lambda$ has dimension
$\binom{D+3}{3}$.
\end{lemma}
\begin{proof} By  \eqref{eq:pfSIZE} and Definition \ref{def:Lambda}.
\end{proof}

\begin{lemma} 
\label{lem:EVec}
For a profile $(r,s,t,u) \in \mathcal P_D$ the following {\rm (i)--(iii)} hold:
\begin{enumerate}
\item[\rm (i)] $A^{*(1)} B(r,s,t,u) = \theta^*_{t+u} B(r,s,t,u)$;
\item[\rm (ii)] $A^{*(2)} B(r,s,t,u) = \theta^*_{u+s} B(r,s,t,u)$;
\item[\rm (iii)] $A^{*(3)} B(r,s,t,u) = \theta^*_{s+t} B(r,s,t,u)$.
\end{enumerate}
\end{lemma}
\begin{proof} By Definition \ref{def:dualmaps} and Remark \ref{rem:aside}.
\end{proof}

\noindent Next, we describe a $\Lambda$-basis that is dual to the $\Lambda$-basis in Definition \ref{def:Lambda} with respect to $\langle \,,\,\rangle$.
\begin{definition}
\label{def:Brstu}
\rm
For a profile $(r,s,t,u) \in {\mathcal P}_D$ define a vector
\begin{align*}
\tilde B(r,s,t,u) = \frac{B(r,s,t,u)}{n(r,s,t,u)},
\end{align*}
where $n(r,s,t,u)$ is from Definition \ref{def:norm}.
\end{definition}

\begin{lemma} 
\label{lem:dualbasis}
The vectors
\begin{align*}
\tilde B(r,s,t,u) \qquad \qquad (r,s,t,u) \in {\mathcal P}_D              
\end{align*}
form an orthogonal basis for $\Lambda$. 
\end{lemma}
\begin{proof} By the construction and Lemma \ref{lem:Borthog}.
\end{proof}

\begin{lemma} 
\label{lem:dualbasis2}
The $\Lambda$-basis 
\begin{align*}
B(r,s,t,u) \qquad \qquad (r,s,t,u) \in {\mathcal P}_D                 
\end{align*}
and the $\Lambda$-basis
\begin{align*}
\tilde B(r,s,t,u) \qquad \qquad (r,s,t,u) \in {\mathcal P}_D                  
\end{align*}
are dual 
with respect to $\langle \,,\,\rangle$.
\end{lemma}
\begin{proof}  For $(r,s,t,u) \in {\mathcal P}_D$ we have
\begin{align*}
\Bigl \langle \tilde B(r,s,t,u), B(r,s,t,u) \Bigr \rangle = \biggl \langle \frac{B(r,s,t,u)}{n(r,s,t,u)}, B(r,s,t,u) \biggr \rangle = \frac{\Vert B(r,s,t,u) \Vert^2}{n(r,s,t,u)} = 1.
\end{align*}
\end{proof}

\noindent For notational convenience, for $r,s,t,u \in \mathbb Z$ we define $\tilde B(r,s,t,u)=0$ if $(r,s,t,u) \not\in \mathcal P_D$.
\medskip

\noindent In the next result, we describe how the maps in \eqref{eq:gens} act on the $\Lambda$-basis from Lemma \ref{lem:dualbasis}.

\begin{proposition}
 \label{lem:nAactionC} 
 For a profile $(r,s,t,u) \in {\mathcal P}_D$  the following {\rm (i)--(vi)} hold.
\begin{enumerate}
\item[\rm (i)] The vector
\begin{align*}
A^{(1)} \tilde B(r,s,t,u)
\end{align*}
 is a linear combination with the following terms and coefficients:
\begin{align*} 
\begin{tabular}[t]{c|c}
{\rm Term }& {\rm Coefficient} 
 \\
 \hline
   $ \tilde B(r-1,s+1,t,u)$& $C(r,t,s,u) $   \\
 $ \tilde B(r+1,s-1,t,u)$   & $C(s,u,r,t)$\\
  $ \tilde B(r,s,t-1,u+1)$  & $C(t,s,u,r)$ \\
  $ \tilde B(r,s,t+1,u-1) $& $C(u,r,t,s)$
   \end{tabular}
\end{align*}
\item[\rm (ii)] 
the vector
\begin{align*}
A^{(2)} \tilde B(r,s,t,u) 
\end{align*}
 is a linear combination with the following terms and coefficients:
\begin{align*} 
\begin{tabular}[t]{c|c}
{\rm Term }& {\rm Coefficient} 
 \\
 \hline
   $ \tilde B(r-1,s,t+1,u) $& $ C(r,u,t,s)$   \\
 $ \tilde B(r,s-1,t,u+1) $   & $C(s,r,u,t) $\\
  $ \tilde B(r+1,s,t-1,u)$  & $C(t,s,r,u) $ \\
  $ \tilde B(r,s+1,t,u-1) $& $C(u,t,s,r) $
   \end{tabular}
\end{align*}
\item[\rm (iii)] 
the vector
\begin{align*}
A^{(3)} \tilde B(r,s,t,u)
\end{align*}
 is a linear combination with the following terms and coefficients:
\begin{align*} 
\begin{tabular}[t]{c|c}
{\rm Term }& {\rm Coefficient} 
 \\
 \hline
    $ \tilde B(r-1,s,t,u+1)$& $C(r,s,u,t)  $   \\
 $\tilde B(r,s-1,t+1,u)  $   & $C(s,u,t,r)  $\\
    $ \tilde B(r,s+1,t-1,u) $  & $C(t,r,s,u)  $ \\
  $ \tilde B(r+1,s,t,u-1) $& $C(u,t,r,s)  $
   \end{tabular}
\end{align*}
\item[\rm (iv)] $A^{*(1)} \tilde B(r,s,t,u) = \theta^*_{t+u} \tilde B(r,s,t,u)$,
\item[\rm (v)] $A^{*(2)} \tilde B(r,s,t,u) = \theta^*_{u+s} \tilde B(r,s,t,u)$,
\item[\rm (vi)] $A^{*(3)} \tilde B(r,s,t,u) = \theta^*_{s+t} \tilde B(r,s,t,u)$.
\end{enumerate}
\end{proposition}
\begin{proof} (i) By Proposition \ref{thm:2homALT} and Definitions \ref{def:Bcheck}, \ref{def:Brstu}. \\
\noindent (ii), (iii) By (i) and $S_3$-symmetry. \\
\noindent (iv)--(vi) By Lemma \ref{lem:EVec} and
Definition \ref{def:Brstu}.
\end{proof} 

\begin{corollary}
\label{cor:LambdaModule}
The vector space $\Lambda$ is invariant under the maps listed in \eqref{eq:gens}. 
\end{corollary}
\begin{proof} By Lemma \ref{lem:dualbasis} and Proposition \ref{lem:nAactionC}.
\end{proof}

\medskip
\noindent  
By Proposition \ref{lem:nAactionC}(iv)--(vi), the $\Lambda$-basis from Lemma \ref{lem:dualbasis} consists of common eigenvectors for
 $A^{*(1)}$, $A^{*(2)}$, $A^{*(3)}$.
Our next general goal is to obtain a $\Lambda$-basis consisting of common
eigenvectors for $A^{(1)}$, $A^{(2)}$, $A^{(3)}$.

\begin{definition}\rm 
\label{def:elementQhij} 
(See   \cite[Section~4]{norton}.)         For $0 \leq h,i,j\leq D$ define a vector
\begin{align*}
Q_{h,i,j} = \vert X \vert   \sum_{x \in X} E_h x \otimes E_i x \otimes E_j x.
\end{align*}
\end{definition}

\begin{lemma}
\label{lem:Qfacts} 
{\rm   (See   \cite[Lemma~4.2]{norton}.)} 
The  vectors 
\begin{align*}
Q_{h,i,j} \qquad \qquad (0 \leq h,i,j\leq D)
\end{align*}
are mutually orthogonal. 
Moreover, 
\begin{align*}
 \Vert Q_{h,i,j} \Vert^2 = \vert X \vert m_h q^h_{i,j} \qquad \qquad (0 \leq h,i,j\leq D).
\end{align*}
\end{lemma}

\begin{corollary} 
\label{cor:pp} For $0 \leq h,i,j\leq D$ the vector $Q_{h,i,j}\not=0$  if and only if $(h,i,j) \in \mathcal P'_D$.
\end{corollary}
\begin{proof}
By Lemmas \ref{lem:nonzINT}, \ref{lem:Qfacts} and since $p^h_{i,j} = q^h_{i,j}$.
\end{proof}

\begin{definition}\label{lem:BstarVector} \rm For a profile $(r,s,t,u) \in \mathcal P_D$ define a vector
\begin{align*}
B^*(r,s,t,u) = Q_{h,i,j},
\end{align*}
\noindent where
\begin{align}
h = t+u, \qquad \qquad i = u+s, \qquad \qquad j = s+t.       \label{eq:hij}
\end{align}
\end{definition}

\noindent  Recall  from  \eqref{eq:oneBF} the vector ${\bf 1}=\sum_{x \in X} x$. We will be discussing the vector ${\bf 1}^{\otimes 3} = {\bf 1} \otimes {\bf 1} \otimes {\bf 1}$. Note that
\begin{align*}
{\bf 1}^{\otimes 3} = \sum_{x,y,z \in X} x \otimes y \otimes z.
\end{align*}

\begin{example} \label{lem:BsD000}
\rm (See \cite[Lemma~9.18]{S3}.) We have
\begin{align*}
B^*(D,0,0,0) = \vert X \vert^{-1} {\bf 1}^{\otimes 3}.
\end{align*}
\end{example} 

\begin{lemma} \label{lem:BsOrthog} The  vectors 
\begin{align*}
B^*(r,s,t,u) \qquad \qquad (r,s,t,u) \in \mathcal P_D         
\end{align*}
are mutually orthogonal. Moreover
\begin{align*}
\Vert B^*(r,s,t,u) \Vert^2 = n(r,s,t,u) \qquad \qquad (r,s,t,u) \in \mathcal P_D.      
\end{align*}
\end{lemma} 
\begin{proof} This is a reformulation of Lemma \ref{lem:Qfacts}, using Lemma \ref{lem:norm} and Definition \ref{lem:BstarVector} along with
 $q^h_{i,j} = p^h_{i,j}$ $(0 \leq h,i,j\leq D)$ and $m_h = k_h$  $(0 \leq h \leq D)$.
\end{proof}

\begin{proposition} \label{lem:basisBstar}  The  vectors 
\begin{align}
B^*(r,s,t,u) \qquad \qquad (r,s,t,u) \in \mathcal P_D  \label{eq:BsList}
\end{align}
form an orthogonal basis for $\Lambda$.
\end{proposition}
\begin{proof} We first show that $\Lambda$ contains the vectors listed in \eqref{eq:BsList}. For $(r,s,t,u) \in \mathcal P_D$ we show that $B^*(r,s,t,u) \in \Lambda$. 
By Definition \ref{lem:BstarVector} we have $B^*(r,s,t,u) = Q_{h,i,j}$ where $Q_{h,i,j}$ is from Definition \ref{def:elementQhij} and $h,i,j$ are from \eqref{eq:hij}. We saw in Section 3 
that $A$ generates $M$ and for $0 \leq \ell \leq D$ the primitive idempotent
$E_\ell$ is contained in $M$. Therefore, there exists a polynomial $g_\ell$ in one variable such that $E_\ell = g_\ell (A)$.
By Example  \ref{ex:BDzzz} and Definition  \ref{def:elementQhij},
\begin{align*}
B^*(r,s,t,u) = \vert X \vert g_h\bigl (   A^{(1)}   \bigr) g_i\bigl (   A^{(2)}   \bigr)  g_j\bigl (   A^{(3)}   \bigr)  B(D,0,0,0).
\end{align*}
The vector space $\Lambda$ contains $B(D,0,0,0)$ and is invariant under $A^{(1)}$, $A^{(2)}$, $A^{(3)}$. Therefore $B^*(r,s,t,u) \in \Lambda$.
We have shown that $\Lambda$ contains the vectors listed in \eqref{eq:BsList}.
The dimension of $\Lambda$ is $\binom{D+3}{3}$, and this is the number of vectors listed in \eqref{eq:BsList}.
By Lemma \ref{lem:BsOrthog}, the vectors listed in \eqref{eq:BsList}  are nonzero and mutually orthogonal. 
 By these comments,
the vectors listed in \eqref{eq:BsList} form an orthogonal basis for $\Lambda$.
\end{proof}

\begin{lemma} 
\label{lem:dualEV}
For a profile $(r,s,t,u) \in \mathcal P_D$ the following {\rm (i)--(iii)} hold:
\begin{enumerate}
\item[\rm (i)] $A^{(1)} B^*(r,s,t,u) = \theta_{t+u} B^*(r,s,t,u)$;
\item[\rm (ii)] $A^{(2)} B^*(r,s,t,u) = \theta_{u+s} B^*(r,s,t,u)$;
\item[\rm (iii)] $A^{(3)} B^*(r,s,t,u) = \theta_{s+t} B^*(r,s,t,u)$.
\end{enumerate}
\end{lemma}
\begin{proof} (i) Let $h,i,j$ be as in  \eqref{eq:hij}. Using Lemma \ref{lem:clarify} and Definitions \ref{def:elementQhij}, \ref{lem:BstarVector} along with $AE_\ell = \theta_\ell E_\ell$ for $0 \leq \ell \leq D$,
\begin{align*}
A^{(1)} B^*(r,s,t,u) &= A^{(1)} Q_{h,i,j} 
                             = A^{(1)} \vert X \vert \sum_{x \in X} E_h x \otimes E_i x \otimes E_j x \\
                              &=  \vert X \vert \sum_{x \in X} A E_h x \otimes E_i x \otimes E_j x 
                              =  \theta_h \vert X \vert \sum_{x \in X} E_h x \otimes E_i x \otimes E_j x \\
                              &= \theta_h Q_{h,i,j} 
                              = \theta_{t+u} B^*(r,s,t,u).
\end{align*}
\noindent (ii), (iii) Similar to the proof of (i).
\end{proof}

\noindent Next, we describe the $\Lambda$-basis that is dual to the $\Lambda$-basis in Proposition \ref{lem:basisBstar}
with respect to $\langle \,,\,\rangle$.

\begin{definition} \label{lem:Bsdual} \rm For a profile $(r,s,t,u) \in \mathcal P_D$ define a vector
\begin{align*}
\tilde B^*(r,s,t,u) = \frac{B^*(r,s,t,u)}{n(r,s,t,u)}.
\end{align*}
\end{definition}

\begin{lemma} \label{lem:LamBasis} The  vectors
\begin{align*} 
\tilde B^*(r,s,t,u) \qquad \qquad (r,s,t,u) \in \mathcal P_D 
\end{align*}
form an orthogonal  basis for $\Lambda$.
\end{lemma}
\begin{proof} By Proposition \ref{lem:basisBstar} and Definition \ref{lem:Bsdual}.
\end{proof}

\begin{lemma} \label{lem:dualBstar} The $\Lambda$-basis
\begin{align*}
B^*(r,s,t,u) \qquad \qquad (r,s,t,u) \in \mathcal P_D 
\end{align*}
and the 
$\Lambda$-basis
\begin{align*}
\tilde B^*(r,s,t,u) \qquad \qquad (r,s,t,u) \in \mathcal P_D 
\end{align*}
are dual with respect to $\langle \,,\,\rangle$.
\end{lemma}
\begin{proof} By Lemma  \ref{lem:BsOrthog} and Definition \ref{lem:Bsdual}.
\end{proof}

\noindent The following result is a variation on   \cite[Lemma~9.13]{S3}.

\begin{lemma} \label{lem:dualSumBs}  We have
\begin{align}
B^*(D,0,0,0)= \vert X \vert^{-1} \sum_{(r,s,t,u) \in \mathcal P_D} B(r,s,t,u). \label{eq:Norm1}
\end{align}
\end{lemma}
\begin{proof} By Lemma \ref{lem:interp}, Definition \ref{def:profile}, and Example \ref{lem:BsD000},
\begin{align*}
B^*(D,0,0,0) &= \vert X \vert^{-1} {\bf 1}^{\otimes 3}  
= \vert X \vert^{-1} \sum_{x,y,z \in X} x \otimes y \otimes z \\
& = \vert X \vert^{-1}  \sum_{(r,s,t,u) \in \mathcal P_D} B(r,s,t,u).
\end{align*}
\end{proof}

\noindent The following result is a variation on  \cite[Lemma~9.18]{S3}.
\begin{lemma} \label{lem:BstarSUM}  We have
\begin{align}
B(D,0,0,0)= \vert X \vert^{-1} \sum_{(r,s,t,u) \in \mathcal P_D} B^*(r,s,t,u).   \label{eq:Norm2}
\end{align}
\end{lemma}
\begin{proof}
 Using Example  \ref{ex:BDzzz} and $I = \sum_{\ell=0}^D E_\ell $,
\begin{align*}
B(D,0,0,0) &= \sum_{x \in X} x \otimes x \otimes x \\
                 &= \sum_{x \in X} \Biggl ( \sum_{h=0}^D E_h x\Biggr ) \otimes  \Biggl ( \sum_{i=0}^D E_i x\Biggr ) \otimes  \Biggl ( \sum_{j=0}^D E_j x\Biggr ) \\
                 &= \sum_{x \in X} \sum_{h=0}^D \sum_{i=0}^D \sum_{j=0}^D E_h x \otimes E_i x \otimes E_j x \\
                  &=  \sum_{h=0}^D \sum_{i=0}^D \sum_{j=0}^D \sum_{x \in X} E_h x \otimes E_i x \otimes E_j x \\
                  &=  \vert X \vert^{-1}  \sum_{h=0}^D \sum_{i=0}^D \sum_{j=0}^D Q_{h,i,j} \\
                   &=  \vert X \vert^{-1}  \sum_{(h,i,j) \in \mathcal P'_D} Q_{h,i,j} \\
                  & = \vert X \vert^{-1} \sum_{(r,s,t,u) \in \mathcal P_D} B^*(r,s,t,u).
\end{align*}
\end{proof}

\noindent For notational convenience, for $r,s,t,u \in \mathbb Z$ we define $\tilde B^*(r,s,t,u)=0$ if $(r,s,t,u) \not\in \mathcal P_D$.
\medskip

\noindent In the next result, we describe how the maps in \eqref{eq:gens} act on the $\Lambda$-basis given in Lemma  \ref{lem:LamBasis}.

\begin{proposition}
 \label{prop:dualAct}
 For a profile $(r,s,t,u) \in \mathcal P_D$ the following {\rm (i)--(vi)} hold:
\begin{enumerate}
\item[\rm (i)] $A^{(1)} \tilde B^*(r,s,t,u) = \theta_{t+u} \tilde B^*(r,s,t,u)$;
\item[\rm (ii)] $A^{(2)} \tilde B^*(r,s,t,u) = \theta_{u+s} \tilde B^*(r,s,t,u)$;
\item[\rm (iii)] $A^{(3)} \tilde B^*(r,s,t,u)  = \theta_{s+t} \tilde B^*(r,s,t,u)$;
\item[\rm (iv)] the vector
\begin{align*}
A^{*(1)} \tilde B^*(r,s,t,u)
\end{align*}
 is a linear combination with the following terms and coefficients:
\begin{align*} 
\begin{tabular}[t]{c|c}
{\rm Term }& {\rm Coefficient} 
 \\
 \hline
   $ \tilde B^*(r-1,s+1,t,u)$& $C(r,t,s,u)$   \\
 $\tilde B^*(r+1,s-1,t,u)$   & $C(s,u,r,t)$\\
  $ \tilde B^*(r,s,t-1,u+1)$  & $C(t,s,u,r)$ \\
  $ \tilde B^*(r,s,t+1,u-1) $& $C(u,r,t,s)$
   \end{tabular}
\end{align*}
\item[\rm (v)] 
the vector
\begin{align*}
A^{*(2)} \tilde B^*(r,s,t,u) 
\end{align*}
 is a linear combination with the following terms and coefficients:
\begin{align*} 
\begin{tabular}[t]{c|c}
{\rm Term }& {\rm Coefficient} 
 \\
 \hline
   $ \tilde B^*(r-1,s,t+1,u) $& $ C(r,u,t,s)$   \\
 $ \tilde B^*(r,s-1,t,u+1) $   & $C(s,r,u,t)$\\
  $ \tilde B^*(r+1,s,t-1,u)$  & $C(t,s,r,u)$ \\
  $ \tilde B^*(r,s+1,t,u-1) $& $C(u,t,s,r) $
   \end{tabular}
\end{align*}
\item[\rm (vi)] 
the vector
\begin{align*}
A^{*(3)} \tilde B^*(r,s,t,u)
\end{align*}
 is a linear combination with the following terms and coefficients:
\begin{align*} 
\begin{tabular}[t]{c|c}
{\rm Term }& {\rm Coefficient} 
 \\
 \hline
    $ \tilde B^*(r-1,s,t,u+1)$& $C(r,s,u,t) $   \\
 $ \tilde B^*(r,s-1,t+1,u)  $   & $C(s,u,t,r) $\\
    $ \tilde B^*(r,s+1,t-1,u) $  & $C(t,r,s,u)$ \\
  $ \tilde  B^*(r+1,s,t,u-1) $& $C(u,t,r,s) $
   \end{tabular}
\end{align*}
\end{enumerate}
\end{proposition}
\noindent The proof of Proposition  \ref{prop:dualAct} is postponed until the end of Section 10.

\begin{remark}\rm In \cite[Section 6]{scaffold} the triply-regular condition is discussed.
By Proposition  \ref{lem:nAactionC}  and \cite[Theorem~6.1(i)]{scaffold} the graph $\Gamma$ is triply-regular.
By Proposition  \ref{prop:dualAct} and \cite[Theorem~6.1(ii)]{scaffold} the graph $\Gamma$ is dual triply-regular.
\end{remark}


\section{Some relations}

\noindent We continue to discuss the  graph $\Gamma$ from Assumption \ref{ASSUME}.
Recall the vector space $\Lambda$ from Definition \ref{def:Lambda}.
In Corollary \ref{cor:LambdaModule}, we saw that  $\Lambda$ is invariant under  the maps listed in \eqref{eq:gens}. 
In this section, we display some relations satisfied by the maps \eqref{eq:gens} acting on $\Lambda$.
\medskip

\noindent Recall the commutator $\lbrack R,S\rbrack=RS-SR$ and the $q$-commutator $\lbrack R,S \rbrack_q = qRS-q^{-1} SR$.

\begin{proposition}
 \label{lem:Nrels} 
  The following relations hold on $\Lambda$:
\begin{enumerate}
\item[\rm (i)] For distinct $i,j \in \lbrace 1,2,3\rbrace$,
\begin{align*}
\lbrack A^{(i)}, A^{(j)} \rbrack=0, \qquad \qquad \lbrack A^{*(i)}, A^{*(j)} \rbrack=0.
\end{align*}
\item[\rm (ii)] For $i \in \lbrace 1,2,3\rbrace$,
\begin{align*}
\lbrack A^{(i)}, A^{*(i)} \rbrack=0.
\end{align*}
\item[\rm (iii)]
For distinct $i,j \in \lbrace 1,2,3\rbrace$,
\begin{align*}
&A^{(i)2} A^{*(j)}  - (q^2+q^{-2}) A^{(i)}  A^{*(j)} A^{(i)} + A^{*(j)} A^{(i)2}  =-H^2 (q^2-q^{-2})^2 A^{*(j)},\\
&A^{*(j)2} A^{(i)}  - (q^2+q^{-2}) A^{*(j)}  A^{(i)} A^{*(j)} + A^{(i)} A^{*(j)2}  =-H^2(q^2-q^{-2})^2  A^{(i)}.
\end{align*}
\item[\rm (iv)] For mutually distinct $h,i,j \in \lbrace 1,2,3\rbrace$,
\begin{align*}
\lbrack A^{(h)}, \lbrack A^{*(i)}, A^{(j)} \rbrack _q \rbrack_q = \lbrack A^{*(h)}, \lbrack A^{(i)}, A^{*(j)} \rbrack_q \rbrack_q.
\end{align*}
\end{enumerate}
\end{proposition}
\begin{proof} (i), (ii) By Definitions \ref{def:maps}, \ref{def:dualmaps} these relations hold on $V^{\otimes 3}$. Therefore these relations hold on $\Lambda$.
\\
\noindent (iii), (iv) To verify these relations, apply each side to a basis vector $\tilde B(r,s,t,u)$ and evaluate the result using
 Proposition  \ref{lem:nAactionC}. More details are given in the Appendix.
 \end{proof}
 
 \begin{remark}\rm Referring to Proposition  \ref{lem:Nrels}, the relations (iii) hold on $V^{\otimes 3}$; this can be shown using
 the methods of \cite[Section~8]{S3}. It is routine to show that the relations (iv) do not hold on $V^{\otimes 3}$ in general.
 \end{remark} 
 
 \begin{remark}\rm Referring to Proposition \ref{lem:Nrels}, the relations (iii) are a special case of the Askey-Wilson relations; see 
 \cite{vidunas, zhedhidden} for a discussion of general Askey-Wilson relations, and \cite{nomTB} for a discussion of the special case.
 \end{remark}
 
 \noindent In the next result, we use Proposition \ref{lem:Nrels} to show that on $\Lambda$, any one of the six generators 
 \begin{align*}
A^{(1)}, \quad A^{(2)}, \quad A^{(3)}, \quad A^{*(1)}, \quad A^{*(2)}, \quad A^{*(3)}
\end{align*}
can be recovered from the other five.

\begin{lemma} 
\label{lem:recover} 
For mutually distinct $h,i,j \in \lbrace 1,2,3 \rbrace$ the following relations hold on $\Lambda$:
\begin{align} 
\label{eq:recover1}
&H^4 (q^2-q^{-2})^4 A^{(h)} = 
\lbrack A^{*(i)} , \lbrack A^{(j)} , \lbrack A^{*(h)}, \lbrack A^{(i)}, A^{*(j)} \rbrack_q \rbrack_q \rbrack_q \rbrack_q, \\
\label{eq:recover2}
&H^4 (q^2-q^{-2})^4 A^{*(h)} = 
\lbrack A^{(i)} , \lbrack A^{*(j)} , \lbrack A^{(h)}, \lbrack A^{*(i)}, A^{(j)} \rbrack_q \rbrack_q \rbrack_q \rbrack_q.
\end{align}
\end{lemma}
 \begin{proof} We first verify \eqref{eq:recover1}. By Proposition \ref{lem:Nrels}(iii), 
 \begin{align}
 \lbrack A^{*(i)}, \lbrack A^{(h)}, A^{*(i)} \rbrack_q \rbrack_q &= H^2(q^2-q^{-2})^2 A^{(h)}, \label{eq:recA} \\
  \lbrack A^{(j)}, \lbrack A^{*(i)}, A^{(j)} \rbrack_q \rbrack_q &= H^2(q^2-q^{-2})^2 A^{*(i)}. \label{eq:recB}
 \end{align}
We may now argue
\begin{align*}
H^4 (q^2-q^{-2})^4 A^{(h)} 
&= H^2 (q^2-q^{-2})^2  \lbrack A^{*(i)}, \lbrack A^{(h)}, A^{*(i)} \rbrack_q \rbrack_q  \qquad \qquad \,\;\hbox{\rm by \eqref{eq:recA}}
\\
&= \lbrack A^{*(i)} , \lbrack A^{(h)} , \lbrack A^{(j)}, \lbrack A^{*(i)}, A^{(j)} \rbrack_q \rbrack_q \rbrack_q \rbrack_q \qquad \qquad \; \hbox{\rm by \eqref{eq:recB}} \\
&= \lbrack A^{*(i)} , \lbrack A^{(j)} , \lbrack A^{(h)}, \lbrack A^{*(i)}, A^{(j)} \rbrack_q \rbrack_q \rbrack_q \rbrack_q \qquad \qquad \;\hbox{since $A^{(h)}, A^{(j)}$ commute} \\
&= \lbrack A^{*(i)} , \lbrack A^{(j)} , \lbrack A^{*(h)}, \lbrack A^{(i)}, A^{*(j)} \rbrack_q \rbrack_q \rbrack_q \rbrack_q \qquad \qquad \hbox{by Proposition \ref{lem:Nrels}(iv)}.
\end{align*}
 We have verified \eqref{eq:recover1}. The verification of \eqref{eq:recover2} is similar.
 \end{proof}

 \noindent In the next section, we will explain what Proposition   \ref{lem:Nrels} and Lemma \ref{lem:recover} have to do with the nonstandard quantum group
  $U'_q(\mathfrak{so}_6)$.

 \section{The nonstandard quantum group $U'_q(\mathfrak{so}_n)$}
 
 In this section,  we fix an integer $n\geq 3$ and do the following.
 First, we motivate things by describing the Lie algebra $\mathfrak{so}_n$. Next, we describe the nonstandard quantum group  $U'_q(\mathfrak{so}_n)$ introduced by Gavrilik and
 Klimyk \cite{gav2}.
 Next, we describe what Proposition   \ref{lem:Nrels} and Lemma \ref{lem:recover} have to do with
  $U'_q(\mathfrak{so}_6)$.
\medskip

\noindent let ${\rm Mat}_n(\mathbb C)$ denote the algebra of $n \times n$ matrices that have all entries in $\mathbb C$.
For $1 \leq i,j\leq n$ define $E_{i,j} \in {\rm Mat}_n(\mathbb C)$ that has $(i,j)$-entry 1 and all other entries $0$.
\medskip

\noindent
The Lie algebra $\mathfrak{gl}_n=\mathfrak{gl}_n(\mathbb C)$ consists of the vector space ${\rm Mat}_n(\mathbb C)$ together with the
Lie bracket $\lbrack R,S \rbrack = RS-SR$. The elements
\begin{align*}
E_{i,j} \qquad \qquad (1 \leq i,j\leq n)
\end{align*}
form a basis for  $\mathfrak{gl}_n$. The dimension of  $\mathfrak{gl}_n$ is $n^2$.
\medskip

\noindent 
For $R \in \mathfrak{gl}_n$ consider the transpose $R^t$. We say that $R$ is {\it antisymmetric} whenever $R^t = -R$.
For $R,S \in  \mathfrak{gl}_n$, if each of $R,S$ is antisymmetric then so is $\lbrack R,S\rbrack$. 
Let  $\mathfrak{so}_n=\mathfrak{so}_n(\mathbb C)$ denote the Lie subalgebra of  $\mathfrak{gl}_n$ consisting of the antisymmetric matrices.
The dimension of  $\mathfrak{so}_n$ is $\binom{n}{2}$.
\medskip

\noindent  Recall from \cite[Theorem~7.36]{carter} the simple Lie algebras over $\mathbb C$ that have finite dimension at least 2:
 \begin{align*}
 A_\ell \;(\ell \geq 1), \quad B_\ell \;(\ell \geq 2), \quad C_\ell \;(\ell \geq 3), \quad D_\ell \;(\ell\geq 4), \quad E_\ell \;(\ell=6,7,8), \quad F_4, \quad G_2.
 \end{align*}
 \begin{lemma} 
 \label{lem:ADE} 
 {\rm (See \cite[Section~21.2]{fulton}.)}
 We give some isomorphisms involving $\mathfrak{so}_n$:
 \begin{enumerate}
 \item[\rm (i)]
 $\mathfrak{so}_3$ is isomorphic to $A_1$;
 \item[\rm (ii)]
  $\mathfrak{so}_4$ is isomorphic to $A_1 \oplus A_1$;
  \item[\rm (iii)]
  $\mathfrak{so}_5$ is isomorphic to $B_2$;
  \item[\rm (iv)]
   $\mathfrak{so}_6$ is isomorphic to $A_3$;
   \item[\rm (v)]
  For odd $n=2r+1\geq 7$,  $\mathfrak{so}_n$ is isomorphic to  $B_r$;
  \item[\rm (vi)]
  For even $n=2r\geq 8$,  $\mathfrak{so}_n$ is isomorphic to $D_r$.
  \end{enumerate}
 \end{lemma}

\noindent Our next goal is to describe a basis for $\mathfrak{so}_n$.

\begin{definition}
\label{def:xijBasis}
\rm
For distinct $i,j\in \lbrace 1,2,\ldots, n\rbrace$ define
\begin{align*}
I_{i,j} = E_{i,j} - E_{j,i}.
\end{align*}
\end{definition}

 \begin{lemma}
 \label{lem:xijRel}
 The elements
 \begin{align*}
 I_{i,j} \qquad \qquad (1 \leq i<j\leq n)
 \end{align*}
 form a basis for  $\mathfrak{so}_n$.
 Moreover, the following relations hold.
       \begin{enumerate}
     \item[\rm (i)] For distinct $i,j \in \lbrace 1,2,\ldots n\rbrace$,
     \begin{align*}
     I_{i,j} = - I_{j,i}.
     \end{align*}
     \item[\rm (ii)] For mutually distinct   $h, i,j \in \lbrace 1,2,\ldots n\rbrace$,
     \begin{align*}
     \lbrack I_{h,i}, I_{i,j} \rbrack = -I_{j,h}.
     \end{align*}
         \item[\rm (iii)] For mutually distinct  $h,i,j,k \in \lbrace 1,2,\ldots n\rbrace$,
     \begin{align*}
     \lbrack I_{h,i}, I_{j,k} \rbrack = 0.
     \end{align*}
     \end{enumerate}
 \end{lemma}
 \begin{proof} This is routinely checked.
 \end{proof}
 
 \noindent We just gave a basis for $\mathfrak{so}_n$.
 Using Lemma  \ref{lem:xijRel}(ii) we may express certain basis elements in terms of others. Applying this idea we find that the
  Lie algebra  $\mathfrak{so}_n$ is generated by $I_{1,2},I_{2,3},\ldots, I_{n-1,n}$. 
 For this generating set, we now give the corresponding presentation of  $\mathfrak{so}_n$ by generators and relations.
 
 \begin{definition}
 \label{def:LL}
 \rm
 Define a Lie algebra $\mathbb L_n$ by generators
 \begin{align*}
 B_i \qquad \qquad (1 \leq i \leq n-1)
 \end{align*}
 and the following relations.
 \begin{enumerate}
 \item[\rm (i)] For $1 \leq i,j\leq n-1$ with $\vert i - j \vert =1$,
 \begin{align*}
 \lbrack B_i, \lbrack B_i, B_j \rbrack \rbrack = - B_j, \qquad \qquad   \lbrack B_j, \lbrack B_j, B_i \rbrack \rbrack = - B_i.
 \end{align*}
 \item[\rm (ii)]  For $1 \leq i,j\leq n-1$ with $\vert i - j \vert \geq 2$,
 \begin{align*}
 \lbrack B_i, B_j \rbrack  = 0.
 \end{align*}
 \end{enumerate}
 \end{definition}
 
 \begin{lemma} 
 \label{lem:solso}
 {\rm (See \cite[Theorem~1]{klimykNS}.)}
 There exists a Lie algebra isomorphism  $\mathbb L_n \to \mathfrak{so}_n$ that sends
 $B_i \mapsto I_{i,i+1}$ for $1 \leq i \leq n-1$.
 \end{lemma}

\noindent Motivated by Definition \ref{def:LL} and Lemma \ref{lem:solso},  we now define $U'_q(\mathfrak{so}_n)$.

\begin{definition}
\label{def:Uqson}
\rm (See  \cite[Section~2]{gav2}.) Assume that $0 \not=q \in \mathbb C$ is not a root of unity.
Define the algebra  $U'_q(\mathfrak{so}_n)$ by generators
 \begin{align*}
 B_i \qquad \qquad (1 \leq i \leq n-1)
 \end{align*}
 and the following relations.
 \begin{enumerate}
 \item[\rm (i)] For $1 \leq i,j\leq n-1$ with $\vert i - j \vert =1$,
 \begin{align*}
 B_i^2 B_j - (q^2+q^{-2}) B_i B_j B_i + B_j B_i^2 &= - B_j, \\
 B_j^2 B_i - (q^2+q^{-2}) B_j B_i B_j+ B_i B_j^2 &= - B_i.
 \end{align*}
 \item[\rm (ii)]  For $1 \leq i,j\leq n-1$ with $\vert i - j \vert \geq 2$,
 \begin{align*}
 \lbrack B_i, B_j \rbrack  = 0.
 \end{align*}
 \end{enumerate}
 \end{definition}
 
 \begin{remark}\rm The notation in \cite[Section~2]{gav2} is different from ours. The scalar $q$ in \cite[Section~2]{gav2} is the same as our $q^2$.
 \end{remark}
 
 \noindent Next, we recall the concept of a PBW basis.
 
 \begin{definition}
\label{def:PBW}
\rm
Let $\mathcal A$ denote an algebra. A {\it Poincar{\'e}-Birkoff-Witt basis} (or {\it PBW basis}) for $\mathcal A$
is a subset $\Omega$ of $\mathcal A$ together with a linear order $\leq $ on $\Omega$ such that the following
is a linear basis for the vector space $\mathcal A$:
\begin{align*}
a_1 a_2 \cdots a_r \qquad \qquad r \in \mathbb N, \qquad a_1, a_2, \ldots, a_r \in \Omega, \qquad a_1 \leq a_2 \leq \cdots \leq a_r.
\end{align*}
\end{definition}
 
  \noindent Our next goal is to describe a PBW basis for  $U'_q(\mathfrak{so}_n)$ that is analogous to the basis for $\mathfrak{so}_n$ given in Lemma \ref{lem:xijRel}.

 \begin{definition}
 \label{def:PBWU}
 \rm (See \cite[Section~2]{klimyk}.)
 For distinct $i,j \in \lbrace 1,2,\ldots, n\rbrace$ we define $I_{i,j} \in  U'_q(\mathfrak{so}_n)$ as follows.
 \begin{enumerate}
 \item[\rm (i)] For $j=i+1$,
 \begin{align*}
 I_{i,i+1} = B_i.
 \end{align*}
  \item[\rm (ii)] For $j\geq i+2$,
 \begin{align*}
 I_{i,j} = \lbrack B_i, I_{i+1,j} \rbrack_q.
 \end{align*}
 \item[\rm (iii)] For $j<i$,
 \begin{align*}
 I_{i,j} = -I_{j,i}.
 \end{align*}
 \end{enumerate}
 \end{definition}
 
 \begin{lemma}
 \label{lem:PBWdetails}
 {\rm (See \cite[Section~2]{klimyk}.)}
 A PBW basis for  $U'_q(\mathfrak{so}_n)$ is obtained by the elements 
 \begin{align*}
 I_{i,j} \qquad \qquad (1 \leq i<j\leq n)
 \end{align*}
 in the following order:
 \begin{align*}
 I_{1,2}<I_{1,3} < \cdots < I_{1,n} <
 I_{2,3}< I_{2,4} < \cdots < I_{2,n} <
 \cdots <I_{n-1,n}.
 \end{align*}
 \end{lemma}
 
 \noindent Next we describe some relations satisfied by the elements $I_{i,j}$ from Definition \ref{def:PBWU}. To facilitate this description, we give a definition.
 \begin{definition}
 \label{def:clockwise} 
 \rm
 Consider a regular $n$-gon $P_n$ with vertices labeled clockwise $1,2,\ldots, n$. We orient the edges $1 \to 2 \to 3 \to \cdots \to n \to 1$. 
  Consider a sequence of at least 3 mutually distinct vertices of $P_n$, written $v_1, v_2, \ldots, v_t$ $(3 \leq t \leq n)$.
 Let $p$ denote the directed path of length $n-1$ that starts at $v_1$ and runs clockwise around $P_n$.
 The sequence $v_1, v_2, \ldots, v_t$ is said to  {\it run clockwise} whenever the path $p$ encounters $v_1,v_2, \ldots, v_t$ in that order.
 The sequence  $v_1, v_2, \ldots, v_t$ is said to run {\it counter-clockwise} (or {\it c-clockwise}) whenever the inverted sequence $v_t, \ldots, v_2, v_1$ is runs clockwise.
 For distinct vertices $i,j$ in $P_n$, by the {\it diagonal $\overline{ij}$} we mean the line segment with endpoints $i,j$.
 \end{definition}
 
 \begin{lemma}
 \label{lem:Urelations}
 {\rm (See \cite[Section~2]{klimyk}.)}
 The following relations are satisfied by the elements $I_{i,j} \in U'_q(\mathfrak{so}_n)$ from Definition  \ref{def:PBWU}.
   \begin{enumerate}
     \item[\rm (i)] For distinct $i,j \in \lbrace 1,2,\ldots n\rbrace$,
     \begin{align*}
     I_{i,j} = - I_{j,i}.
     \end{align*}
     \item[\rm (ii)] For mutually distinct   $h, i,j \in \lbrace 1,2,\ldots n\rbrace$,
     \begin{align*}
     \lbrack I_{h,i}, I_{i,j} \rbrack_q &= -I_{j,h} \qquad \hbox{if the sequence $h,i,j$ runs clockwise;} \\
      \lbrack I_{h,i}, I_{i,j} \rbrack_{q^{-1}} &= -I_{j,h} \qquad \hbox{if the sequence $h,i,j$ runs c-clockwise.} 
     \end{align*}
         \item[\rm (iii)] For mutually distinct  $h,i,j,k \in \lbrace 1,2,\ldots n\rbrace$, 
              \begin{align*}
     \lbrack I_{h,i}, I_{j,k} \rbrack &= 0 \qquad \hbox{if the diagonals $\overline{hi}$ and $\overline{jk}$ do not overlap};\\
        \lbrack I_{h,i}, I_{j,k} \rbrack &= (q^{-2}-q^{2})(I_{h,j} I_{i,k} + I_{j,i} I_{k,h}) \qquad \hbox{if the sequence $h,j,i,k$ runs clockwise};\\
          \lbrack I_{h,i}, I_{j,k} \rbrack &= (q^2-q^{-2})(I_{h,j} I_{i,k} + I_{j,i} I_{k,h}) \qquad \hbox{if the sequence $h,j,i,k$ runs c-clockwise}.
     \end{align*}
     \end{enumerate}
\end{lemma}

 \begin{remark} \rm There is typo in  \cite[line~(2.8)]{klimyk}.  In that line the left-hand side should be a commutator instead of a $q$-commutator.
 \end{remark}
 
 \noindent In Lemma  \ref{lem:Urelations} we see a $\mathbb Z_n$-cyclic symmetry among the relations. We now make this symmetry more explicit.
 
 \begin{definition}
 \label{def:Bn}
 \rm For notational convenience, define $B_n \in U'_q(\mathfrak{so}_n)$ by
 \begin{align*}
 B_n = I_{n,1},
 \end{align*}
  where $I_{n,1}$ is from Definition  \ref{def:PBWU}.
 \end{definition}
 
 \begin{lemma}
 \label{lem:BnSym}
 The following relations hold in  $U'_q(\mathfrak{so}_n)$.
 \begin{enumerate}
 \item[\rm (i)] For $i \in \lbrace 1,n-1\rbrace$,
 \begin{align*}
 B^2_i B_n- (q^2+q^{-2}) B_i B_n B_i + B_n B^2_i &= - B_n,\\
  B^2_n B_i - (q^2+q^{-2}) B_n B_i B_n + B_i B^2_n &= - B_i.
\end{align*}
\item[\rm (ii)] For $2 \leq i \leq n-2$,
\begin{align*}
\lbrack B_n, B_i \rbrack=0.
\end{align*}
\end{enumerate}
 \end{lemma}
 \begin{proof} (i) We first assume that $i=1$.
 By Lemma  \ref{lem:Urelations}(ii) and Definition  \ref{def:Bn},
 \begin{align*}
 -B_n &= - I_{n,1} = \lbrack I_{1,2}, I_{2,n} \rbrack_q = -  \lbrack I_{1,2}, \lbrack I_{n,1}, I_{1,2} \rbrack_q  \rbrack_q 
 =   -  \lbrack B_1, \lbrack B_n, B_1 \rbrack_q  \rbrack_q \\
 &= B^2_1 B_n- (q^2+q^{-2}) B_1 B_n B_1 + B_n B^2_1.
 \end{align*}
 Also by Lemma  \ref{lem:Urelations}(ii) and Definition  \ref{def:Bn},
 \begin{align*}
 -B_1 &= - I_{1,2} = \lbrack I_{2,n}, I_{n,1} \rbrack_q = -   \lbrack  \lbrack I_{n,1}, I_{1,2} \rbrack_q, I_{n,1} \rbrack    \rbrack_q
 = -  \lbrack  \lbrack B_n, B_1 \rbrack_q, B_n \rbrack    \rbrack_q \\
 &= B^2_n B_1- (q^2+q^{-2}) B_n B_1 B_n + B_1 B^2_n.
 \end{align*}
 We have verified the result for $i=1$. The verification for $i=n-1$  is similar.\\
 \noindent (ii) By Lemma  \ref{lem:Urelations}(iii) and Definition  \ref{def:Bn}.
 \end{proof}
 
 \begin{lemma}
 \label{lem:ZSym} There exists an automorphism $\rho$ of  $U'_q(\mathfrak{so}_n)$ that sends
 $B_i \mapsto B_{i+1}$ for $1 \leq i \leq n-1$ and $B_n \mapsto B_1$.
 For distinct $i,j \in \lbrace 1,2,\ldots, n\rbrace$ this automorphism sends $I_{i,j} \mapsto I_{i+1,j+1}$
 where we understand $I_{i,n+1}=I_{i,1}$ and $I_{n+1,j} = I_{1,j}$.
 \end{lemma}
 \begin{proof} For $1 \leq i \leq n-1$ define $\mathcal B_i = B_{i+1}$. By Definition \ref{def:Uqson} and Lemma  \ref{lem:BnSym}, the elements $\lbrace {\mathcal B}_i \rbrace_{i=1}^{n-1}$
 satisfy
 the following (i), (ii).
 \begin{enumerate}
 \item[\rm (i)] For $1 \leq i,j\leq n-1$ with $\vert i - j \vert =1$,
 \begin{align*}
 {\mathcal B}_i^2 {\mathcal B}_j - (q^2+q^{-2}) {\mathcal B}_i {\mathcal B}_j {\mathcal B}_i + {\mathcal B}_j {\mathcal B}_i^2 &= - {\mathcal B}_j, \\
{\mathcal  B}_j^2 {\mathcal B}_i - (q^2+q^{-2}) {\mathcal B}_j {\mathcal B}_i {\mathcal B}_j+ {\mathcal B}_i {\mathcal B}_j^2 &= - {\mathcal B}_i.
 \end{align*}
 \item[\rm (ii)]  For $1 \leq i,j\leq n-1$ with $\vert i - j \vert \geq 2$,
 \begin{align*}
 \lbrack {\mathcal B}_i, {\mathcal B}_j \rbrack  = 0.
 \end{align*}
 \end{enumerate}
Comparing these relations with the relations in Definition \ref{def:Uqson}, we obtain an algebra homomorphism $\rho: U'_q(\mathfrak{so}_n) \to U'_q(\mathfrak{so}_n)$ that sends
$B_i \mapsto {\mathcal B}_i = B_{i+1}$ for $1 \leq i \leq n-1$. We now show that $\rho $ sends $B_n \mapsto B_1$. 
By Definitions  \ref{def:PBWU}, \ref{def:Bn}  we have
\begin{align*}
B_n = I_{n,1} = -  \lbrack B_1, \lbrack B_2, \ldots, \lbrack B_{n-2}, B_{n-1} \rbrack_q \cdots \rbrack_q  \rbrack_q.
\end{align*}
By Lemma  \ref{lem:Urelations} and the construction,
 \begin{align*}
B_1 = I_{1,2} = -  \lbrack B_2, \lbrack B_3, \ldots, \lbrack B_{n-1}, B_{n} \rbrack_q \cdots \rbrack_q  \rbrack_q.
\end{align*}
By these comments, $\rho$ sends $B_n \mapsto B_1$. The map $\rho^n$ fixes $B_i$ for $1 \leq i \leq n$, so $\rho^n = {\rm id}$. 
Consequently  $\rho$ is invertible
and hence a bijection. We have shown that $\rho$ is an automorphism of $U'_q(\mathfrak{so}_n)$ that sends $B_i \mapsto  B_{i+1}$ for $1 \leq i \leq n-1$
and $B_n \mapsto B_1$. The last assertion of the lemma statement is checked using Lemma  \ref{lem:Urelations}.
 \end{proof}
 
\noindent We return our attention to the graph $\Gamma$ from Assumption \ref{ASSUME}. Our next goal is to
explain what Proposition   \ref{lem:Nrels} and Lemma \ref{lem:recover} have to do with $U'_q(\mathfrak{so}_6)$. We will turn 
the vector space $\Lambda$ into a  $U'_q(\mathfrak{so}_6)$-module in two ways.
 \begin{theorem}
 \label{thm:main}
 The vector space $\Lambda$ becomes a $U'_q(\mathfrak{so}_6)$-module on which
 \begin{align*}
& B_1 = \frac{A^{(1)}}{H(q^2-q^{-2})}, \qquad \qquad
 B_3 = \frac{A^{(2)}}{H(q^2-q^{-2})}, \qquad \qquad 
 B_5 = \frac{A^{(3)}}{H(q^2-q^{-2})}, \\
& B_2 = \frac{A^{*(3)}}{H(q^{-2}-q^{2})}, \qquad \qquad
 B_4 = \frac{A^{*(1)}}{H(q^{-2}-q^{2})}, \qquad \qquad
  B_6 = \frac{A^{*(2)}}{H(q^{-2}-q^{2})}.
 \end{align*}
 \end{theorem}
 \begin{proof} Define
  \begin{align*}
&{\sf B}_1 = \frac{A^{(1)}}{H(q^2-q^{-2})}, \qquad \qquad
{\sf B}_3 = \frac{A^{(2)}}{H(q^2-q^{-2})}, \qquad \qquad 
{\sf B}_5 = \frac{A^{(3)}}{H(q^2-q^{-2})}, \\
&{\sf  B}_2 = \frac{A^{*(3)}}{H(q^{-2}-q^{2})}, \qquad \qquad
{\sf B}_4 = \frac{A^{*(1)}}{H(q^{-2}-q^{2})}, \qquad \qquad
{\sf  B}_6 = \frac{A^{*(2)}}{H(q^{-2}-q^{2})}.
 \end{align*}
  By Proposition \ref{lem:Nrels}, on $\Lambda$
  the elements $\lbrace {\sf B}_i \rbrace_{i=1}^{5}$
 satisfy
 the following (i), (ii).
 \begin{enumerate}
 \item[\rm (i)] For $1 \leq i,j\leq 5$ with $\vert i - j \vert =1$,
 \begin{align*}
 {\sf B}_i^2 {\sf B}_j - (q^2+q^{-2}) {\sf B}_i {\sf B}_j {\sf B}_i + {\sf B}_j {\sf B}_i^2 &= - {\sf B}_j, \\
{\sf  B}_j^2 {\sf B}_i - (q^2+q^{-2}) {\sf B}_j {\sf B}_i {\sf B}_j+ {\sf B}_i {\sf B}_j^2 &= - {\sf B}_i.
 \end{align*}
 \item[\rm (ii)]  For $1 \leq i,j\leq 5$ with $\vert i - j \vert \geq 2$,
 \begin{align*}
 \lbrack {\sf B}_i, {\sf B}_j \rbrack  = 0.
 \end{align*}
 \end{enumerate}
Comparing these relations with the relations in Definition \ref{def:Uqson}, we turn $\Lambda$ into a $U'_q(\mathfrak{so}_6)$-module on which
$B_i = {\sf B}_i$ for $1 \leq i \leq 5$. It remains to show that  $B_6={\sf B}_6$ on $\Lambda$.
By Definitions  \ref{def:PBWU}, \ref{def:Bn}  the following holds in  $U'_q(\mathfrak{so}_6)$:
\begin{align} \label{eq:com1}
B_6 = I_{6,1} = - \lbrack B_1, \lbrack B_2, \lbrack B_3, \lbrack B_{4}, B_{5} \rbrack_q  \rbrack_q \rbrack_q \rbrack_q.
\end{align}
By \eqref{eq:recover2}  (with    $h=2$, $i=1$, $j=3$) we see that  on $\Lambda$,
 \begin{align} \label{eq:com2}
{\sf B}_6  = -  \lbrack {\sf B}_1, \lbrack {\sf B}_2, \lbrack {\sf B}_3, \lbrack {\sf B}_{4}, {\sf B}_{5} \rbrack_q \rbrack_q \rbrack_q  \rbrack_q.
\end{align}
Comparing \eqref{eq:com1}, \eqref{eq:com2} we obtain $B_6 ={\sf B}_6$ on $\Lambda$. The result follows.
 \end{proof}

 \begin{lemma} 
 \label{lem:LambdaIRRED}
 The  $U'_q(\mathfrak{so}_6)$-module $\Lambda$ from Theorem \ref{thm:main} is
 irreducible.
 \end{lemma}
 \begin{proof} Let $\mathbb T$ denote the subalgebra of ${\rm End}(V^{\otimes 3})$ generated by the maps listed in \eqref{eq:gens}.
 By Corollary \ref{cor:LambdaModule}, the vector space $\Lambda$ is a $\mathbb T$-submodule of $V^{\otimes 3}$. It suffices to show that the
 $\mathbb T$-module $\Lambda$ is irreducible.
 By \cite[Definition~9.7]{S3}, there exists a unique irreducible  $\mathbb T$-submodule of $V^{\otimes 3}$ that contains
 ${\bf 1}^{\otimes 3}$; this $\mathbb T$-module is called fundamental. By Lemmas \cite[Lemma~9.10]{S3} and \cite[Lemma~9.15]{S3}, the
 fundamental $\mathbb T$-module contains $P_{h,i,j}$ and $Q_{h,i,j}$ for $0 \leq h,i,j\leq D$. In other words, 
the fundamental $\mathbb T$-module contains
 $B(r,s,t,u)$ and $B^*(r,s,t,u)$ for $(r,s,t,u) \in \mathcal P_D$. These vectors span $\Lambda$, so the fundamental $\mathbb T$-module contains $\Lambda$ as a submodule.
 The fundamental $\mathbb T$-module is irreducible, so it is equal  to $\Lambda$.
 We have shown that the $\mathbb T$-module $\Lambda$ is irreducible. The result follows.
 \end{proof}
 
  \begin{theorem}
 \label{thm:main2}
 The vector space $\Lambda$ becomes a $U'_q(\mathfrak{so}_6)$-module on which
 \begin{align*}
& B_1 = \frac{A^{*(1)}}{H(q^2-q^{-2})}, \qquad \qquad
 B_3 = \frac{A^{*(2)}}{H(q^2-q^{-2})}, \qquad \qquad 
 B_5 = \frac{A^{*(3)}}{H(q^2-q^{-2})}, \\
& B_2 = \frac{A^{(3)}}{H(q^{-2}-q^{2})}, \qquad \qquad
 B_4 = \frac{A^{(1)}}{H(q^{-2}-q^{2})}, \qquad \qquad
  B_6 = \frac{A^{(2)}}{H(q^{-2}-q^{2})}.
 \end{align*}
 \end{theorem}
\begin{proof} Similar to the proof of Theorem \ref{thm:main}.
\end{proof}

  \begin{lemma} 
  \label{lem:LambdaIRRED2}
  The  $U'_q(\mathfrak{so}_6)$-module $\Lambda$ from Theorem \ref{thm:main2} is
 irreducible.
 \end{lemma}
 \begin{proof} By Lemma \ref{lem:LambdaIRRED} and the construction.
 \end{proof}

\noindent Shortly we will show that the   $U'_q(\mathfrak{so}_6)$-modules in Theorems \ref{thm:main}, \ref{thm:main2} are isomorphic.
\medskip

\noindent 
We have some comments about the representation theory of $U'_q(\mathfrak{so}_n)$. 
 By \cite[Proposition~5.1]{klimyk}, on each finite-dimensional $U'_q(\mathfrak{so}_n)$-module the generators $\lbrace B_i \vert 1 \leq i \leq n-1, \;i \;\hbox{\rm odd}\rbrace$
 are simultaneously diagonalizable.
  By \cite[Corollary~9.4]{klimyk} each finite-dimensional  $U'_q(\mathfrak{so}_n)$-module is completely reducible;
 this means that the module is a direct sum of irreducible  $U'_q(\mathfrak{so}_n)$-submodules.
 In \cite[Theorem~9.3]{klimyk} the finite-dimensional irreducible $U'_q(\mathfrak{so}_n)$-modules are classified up to isomorphism. 
According to the classification, there are two types of finite-dimensional irreducible $U'_q(\mathfrak{so}_n)$-modules, called classical type and nonclassical type. The type is determined by the form of the
eigenvalues for the generators  $\lbrace B_i \vert 1 \leq i \leq n-1, \;i \;\hbox{\rm odd}\rbrace$  acting on the module. For example, the  $U'_q(\mathfrak{so}_6)$-module $\Lambda$ from Theorem \ref{thm:main} 
or Theorem \ref{thm:main2} 
has classical type; this is verified by comparing \eqref{eq:thsi}  with \cite[Proposition~5.1]{klimyk}. The finite-dimensional irreducible $U'_q(\mathfrak{so}_n)$-modules of classical type are described in \cite[Section~3]{klimyk}. We give some details under the assumption $n=6$.
 The isomorphism classes of  finite-dimensional irreducible $U'_q(\mathfrak{so}_6)$-modules of classical type are in bijection with the 
 3-tuples $(n_1, n_2, n_3)$ such that (i) $2 n_i \in \mathbb Z$ for $i \in \lbrace 1,2,3\rbrace$;
 (ii) $n_i - n_j \in \mathbb Z$ for $i,j\in \lbrace 1,2,3\rbrace$; (iii) $n_1 \geq n_2 \geq \vert n_3\vert$.
 Given a 3-tuple  $(n_1, n_2, n_3)$ that satisfies (i)--(iii) above, the corresponding 
 finite-dimensional irreducible $U'_q(\mathfrak{so}_6)$-module of classical type is constructed in \cite[Section~3]{klimyk}
by giving 
a Gelfand-Tsetlin basis for the module and the action of the generators $\lbrace B_i \rbrace_{i=1}^5$ on the basis.
For a finite-dimensional irreducible $U'_q(\mathfrak{so}_6)$-module of classical type, the corresponding 3-tuple $(n_1,n_2,n_3)$ is
called the {\it highest weight} of the module.
The  $U'_q(\mathfrak{so}_6)$-module $\Lambda$ from Theorem \ref{thm:main} (resp. Theorem \ref{thm:main2}) has highest weight  $(D/2, D/2, D/2)$;
this is verified using Lemma \ref{lem:dualEV}  (resp. Lemma \ref{lem:EVec}) and the description in \cite[Section~1.1]{wenzl}.

\begin{lemma} 
\label{lem:LambdaISO}
The  $U'_q(\mathfrak{so}_6)$-modules  from Theorem \ref{thm:main} and Theorem \ref{thm:main2} are isomorphic.
\end{lemma}
\begin{proof}  We mentioned above the lemma statement that both of these  $U'_q(\mathfrak{so}_6)$-modules have classical type and highest weight $(D/2,D/2,D/2)$.
Since these  $U'_q(\mathfrak{so}_6)$-modules  have the same type and same highest weight, they must be isomorphic by  \cite[Section~3]{klimyk}.
\end{proof}
 
 \noindent {\it Proof of Proposition  \ref{prop:dualAct}.} (i)--(iii) By Lemma \ref{lem:dualEV} and Definition \ref{lem:Bsdual}. \\
 \noindent (iv)--(vi)
 By Lemma \ref{lem:LambdaISO}
 there exists a  $U'_q(\mathfrak{so}_6)$-module isomorphism $K$ from the 
 $U'_q(\mathfrak{so}_6)$-module  in Theorem \ref{thm:main} to the
  $U'_q(\mathfrak{so}_6)$-module in
  Theorem \ref{thm:main2}. By construction $K$ is a $\mathbb C$-linear bijection $\Lambda \to \Lambda$. By Theorems \ref{thm:main} and \ref{thm:main2} the following hold
  on $\Lambda$:
  \begin{align}
  K A^{(i)} = A^{*(i)} K, \qquad \qquad K A^{*(i)} = A^{(i)} K \qquad \qquad i \in \lbrace 1,2,3\rbrace.
  \label{eq:Kconj}
  \end{align}
  By \eqref{eq:Kconj} the map $K^2$ commutes with $A^{(i)}$ and $A^{*(i)}$ for $i \in \lbrace 1,2,3\rbrace$. Therefore  $K^2 \in {\rm Span}(I)$ in view of Lemma \ref{lem:LambdaIRRED}.
   Multiplying $K$ by
  a nonzero complex scalar if necessary, we may assume that $K^2 = I$.
 Let $(r,s,t,u) \in \mathcal P_D$. By Lemma \ref{lem:EVec}, the vector $ B(r,s,t,u)$ is a common eigenvector for $A^{*(1)}$, $A^{*(2)}$,  $A^{*(3)}$ with
 eigenvalues $\theta^*_{t+u}$, $\theta^*_{u+s}$, $\theta^*_{s+t}$ respectively. By this and $\theta_\ell =\theta^*_\ell$ $(0 \leq \ell \leq D)$,
 the vector $K  B(r,s,t,u)$ is a common
 eigenvector for  $A^{(1)}$, $A^{(2)}$, $A^{(3)}$ with
 eigenvalues $\theta_{t+u}$, $\theta_{u+s}$, $\theta_{s+t}$ respectively. By this and Lemma \ref{lem:dualEV}, there exists a nonzero 
$\alpha(r,s,t,u) \in \mathbb C$ such that
\begin{align*}
K B(r,s,t,u) = \alpha(r,s,t,u) B^*(r,s,t,u).
 \end{align*}
We apply $K$ to each side of  \eqref{eq:Norm1} and evaluate the result using \eqref{eq:Norm2}; this yields
\begin{align*}
 \alpha(r,s,t,u)= \frac{1}{\alpha(D,0,0,0)}\qquad \qquad (r,s,t,u) \in \mathcal P_D.
\end{align*}
Setting $(r,s,t,u) =(D,0,0,0)$ we obtain $\alpha(D,0,0,0)^2 = 1$. Replacing $K$ by $-K$ if necessary, we may assume that $\alpha(D,0,0,0)=1$.
Consequently $\alpha(r,s,t,u)=1$ for $(r,s,t,u) \in \mathcal P_D$.
 We have
 \begin{align*} 
 K  B(r,s,t,u) =  B^*(r,s,t,u)     \qquad \qquad (r,s,t,u) \in \mathcal P_D.
 \end{align*}
 By this and Definitions \ref{def:Brstu}, \ref{lem:Bsdual} we obtain
  \begin{align} 
 \label{eq:Kaction}
 K  \tilde B(r,s,t,u) =  \tilde B^*(r,s,t,u)     \qquad \qquad (r,s,t,u) \in \mathcal P_D.
 \end{align}
 To finish the proof, in Proposition  \ref{lem:nAactionC}(i)--(iii) apply $K$ to every vector in the given linear dependency, 
 and evaluate the results using \eqref{eq:Kconj}, \eqref{eq:Kaction}.
\hfill $\Box$ 

\section{Comments}
 
 \noindent In the previous sections, we considered a 2-homogeneous bipartite distance-regular graph $\Gamma$ with diameter $D\geq 3$.
 We assumed that $\Gamma$ is not a hypercube nor a cycle.
 We considered a $Q$-polynomial structure on $\Gamma$.
 We described the corresponding eigenvalue sequence and dual eigenvalue sequence using a nonzero $q \in \mathbb C$ that is not a root of unity.
 Using the standard module $V$ of $\Gamma$, we described a subspace $\Lambda$ of $V^{\otimes 3}$ that has dimension $\binom{D+3}{3}$. We showed how $\Lambda$ becomes an
 irreducible $U'_q(\mathfrak{so}_6)$-module with classical type and highest weight $(D/2,D/2,D/2)$. According to  \cite[Theorem~1.2]{2hom2}
  the graph $\Gamma$ only exists for certain values of $D$ and $q$. Nevertheless, the insight gained from $\Gamma$
 suggests that the following algebraic result holds without restriction on $D$ and $q$.
 
 \begin{proposition}
 \label{thm:finalMain}
 Pick an integer $D\geq 1$. Pick $0 \not=q \in \mathbb C$ that is not a root of unity. Pick any $0 \not=H \in \mathbb C$.
 Define
 the complex scalars $\lbrace \theta_i \rbrace_{i=0}^D$, $\lbrace \theta^*_i \rbrace_{i=0}^D$ as in \eqref{eq:thsi}.
 Define the complex scalars
 \begin{align*}
 C(r,s,t,u) \qquad \qquad (r,s,t,u) \in \mathcal P_D
 \end{align*}
 as in  \eqref{eq:C}. 
 Let $\mathbb V$ denote the finite-dimensional irreducible $U'_q(\mathfrak{so}_6)$-module with classical type and highest weight $(D/2,D/2,D/2)$.
 Define the maps
 \begin{align}
 \label{eq:sixMaps}
 A^{(1)}, \quad   A^{(2)}, \quad  A^{(3)}, \quad 
 A^{*(1)}, \quad  A^{*(2)}, \quad  A^{*(3)}
 \end{align}
 in 
 $ {\rm End}(\mathbb V)$ that satisfy the equations in Theorem  \ref{thm:main} or Theorem  \ref{thm:main2}.
 Then:
 \begin{enumerate}
 \item[\rm (i)] the $U'_q(\mathfrak{so}_6)$-module $\mathbb V$ has a basis
              \begin{align}
              \tilde B(r,s,t,u) \qquad \qquad (r,s,t,u) \in \mathcal P_D         \label{eq:finBasis1}
              \end{align}
              on which the maps \eqref{eq:sixMaps} act according to Proposition  \ref{lem:nAactionC};
 \item[\rm (ii)] the $U'_q(\mathfrak{so}_6)$-module $\mathbb V$ has a basis
              \begin{align}
              \tilde B^*(r,s,t,u) \qquad \qquad (r,s,t,u) \in \mathcal P_D           \label{eq:finBasis2}
              \end{align}
              on which the maps \eqref{eq:sixMaps} act according to Proposition  \ref{prop:dualAct};
 \item[\rm (iii)]  the maps  \eqref{eq:sixMaps} satisfy the relations in Proposition \ref{lem:Nrels};
  \item[\rm (iv)]  the maps  \eqref{eq:sixMaps} satisfy the relations in Lemma \ref{lem:recover}.
 \end{enumerate}
 \end{proposition}
 \begin{proof} (Sketch) We will work with Theorem  \ref{thm:main}; the case of Theorem  \ref{thm:main2} is similar.
 Consider a vector space $\sf V$ of dimension $\binom{D+3}{3}$. Endow $\sf V$ with a basis
 denoted
  \begin{align}
              \tilde {\sf B}(r,s,t,u) \qquad \qquad (r,s,t,u) \in \mathcal P_D.        \label{eq:need1}    
              \end{align}
 Define some maps
  \begin{align}
 {\sf A}^{(1)}, \quad   {\sf A}^{(2)}, \quad  {\sf A}^{(3)}, \quad 
{\sf A}^{*(1)}, \quad  {\sf A}^{*(2)}, \quad  {\sf A}^{*(3)}       \label{eq:sixMapsN}
 \end{align}
 in 
 $ {\rm End}({\sf V})$ that act on the basis vectors \eqref{eq:need1} according to Proposition  \ref{lem:nAactionC}.
One checks that the maps   \eqref{eq:sixMapsN} satisfy the relations 
in Proposition \ref{lem:Nrels} (see the Appendix for details)
  and the relations in Lemma \ref{lem:recover}. By these relations 
  $\sf V$ becomes a $U'_q(\mathfrak{so}_6)$-module that meets the requirements of Theorem \ref{thm:main}.
  One checks that the $U'_q(\mathfrak{so}_6)$-module $\sf V$ is irreducible, with classical type and
  highest weight $(D/2, D/2, D/2)$. Thus the $U'_q(\mathfrak{so}_6)$-modules $\sf V$ and $\mathbb V$ are isomorphic.
  We have shown that the $U'_q(\mathfrak{so}_6)$-module $\mathbb V$    satisfies       (i), (iii), (iv).
  A similar argument shows that the $U'_q(\mathfrak{so}_6)$-module $\mathbb V$    satisfies       (ii), (iii), (iv).
  The result follows.
 \end{proof}
 
 \noindent The following problem is open.
 
 \begin{problem}\rm Referring to the $U'_q(\mathfrak{so}_6)$-module $\mathbb V$ in Proposition  \ref{thm:finalMain},
 find the transition matrices between the $\mathbb V$-basis \eqref{eq:finBasis1} and the $\mathbb V$-basis \eqref{eq:finBasis2}.
 \end{problem}
 
 \section{Appendix}
 
 In this Appendix, we give some details that are used in the proofs of Propositions  \ref{lem:Nrels}, \ref{thm:finalMain}.
 \medskip
 
 \noindent Throughout this Appendix the following assumptions hold. Fix an integer $D\geq 1$. Pick $0 \not=q \in \mathbb C$
 that is not a root of unity. Pick any $0 \not=H \in \mathbb C$.  
  For $i \in \mathbb Z$ define
 the complex scalar  $\theta^*_i$ as in \eqref{eq:thsi}.
Note that
\begin{align}
\theta^*_{i-1} - \beta \theta^*_i + \theta^*_{i+1}  &= 0 \qquad \qquad (i \in \mathbb Z),  \label{eq:rec1}\\
\theta_{i-1}^{*2} - \beta \theta^*_{i-1} \theta^*_i + \theta_i^{*2} &= - H^2 (q^2 - q^{-2})^2 \qquad \qquad (i \in \mathbb Z). \label{eq:rec2}
\end{align}
 For $r,s,t,u \in \mathbb Z$ define the complex scalar
$ C(r,s,t,u)$  as in \eqref{eq:C}. 
 Consider a vector space $\sf V$ of dimension $\binom{D+3}{3}$. Endow $\sf V$ with a basis
 denoted
  \begin{align*}
              \tilde {\sf B}(r,s,t,u) \qquad \qquad (r,s,t,u) \in \mathcal P_D.   
              \end{align*}
 Define some maps
  \begin{align*}
 {\sf A}^{(1)}, \quad   {\sf A}^{(2)}, \quad  {\sf A}^{(3)}, \quad 
{\sf A}^{*(1)}, \quad  {\sf A}^{*(2)}, \quad  {\sf A}^{*(3)}   
 \end{align*}
 in 
 $ {\rm End}({\sf V})$ that act on the above basis according to Proposition  \ref{lem:nAactionC}.
 Our goal is to check that these maps satisfy the relations in Proposition \ref{lem:Nrels}.
 \medskip
 
 \noindent In what follows, let $(r,s,t,u) \in \mathcal P_D$.
\medskip

\noindent The following identities are used to show that ${\sf A}^{(1)} {\sf A}^{(2)} = {\sf A}^{(2)} {\sf A}^{(1)}$ holds at  $ \tilde {\sf B}(r,s,t,u)$:
\begin{align*}
C(r,t,s,u) C(r-1,u,t,s+1)&=C(r,u,t,s)C(r-1,t+1,s,u), \\
C(s,u,r,t)C(t,s-1,r+1,u) &= C(t,s,r,u)C(s,u,r+1,t-1), \\
C(s,u,r,t)C(s-1,r+1,u,t)&= C(s,r,u,t)C(s-1,u+1,r,t), \\
C(r,t,s,u)C(u,t,s+1,r-1)&= C(u,t,s,r)C(r,t,s+1,u-1),\\
C(t,s,u,r)C(t-1,s,r,u+1)&= C(t,s,r,u)C(t-1,s,u,r+1), \\
C(u,r,t,s)C(r,u-1,t+1,s)&= C(r,u,t,s)C(u,r-1,t+1,s),\\
C(u,r,t,s)C(u-1,t+1,s,r)&= C(u,t,s,r)C(u-1,r,t,s+1), \\
C(t,s,u,r)C(s,r,u+1,t-1)&= C(s,r,u,t)C(t,s-1,u+1,r)
\end{align*}
and also
\begin{align*}
&C(r,t,s,u)C(s+1,r-1,u,t)+C(t,s,u,r)C(r,u+1,t-1,s) \\
&= C(r,u,t,s)C(t+1,s,u,r-1)+C(s,r,u,t)C(r,t,s-1,u+1), \\
&C(r,t,s,u)C(t,s+1,r-1,u)+C(t,s,u,r)C(u+1,t-1,s,r) \\
&= C(t,s,r,u)C(r+1,t-1,s,u)+C(u,t,s,r)C(t,s+1,u-1,r), \\
&C(s,u,r,t)C(r+1,u,t,s-1)+C(u,r,t,s)C(s,r,u-1,t+1) \\
&= C(r,u,t,s)C(s,u,r-1,t+1)+C(s,r,u,t)C(u+1,r,t,s-1), \\
&C(s,u,r,t)C(u,t,s-1,r+1)+C(u,r,t,s)C(t+1,s,r,u-1) \\
&= C(t,s,r,u)C(u,r+1,t-1,s)+ C(u,t,s,r)C(s+1,u-1,r,t).
\end{align*}
\noindent The relation
 ${\sf A}^{*(1)} {\sf A}^{*(2)} = {\sf A}^{*(2)} {\sf A}^{*(1)}$ holds by construction.
 \medskip
 
\noindent The relation
 ${\sf A}^{(1)} {\sf A}^{*(1)} = {\sf A}^{*(1)} {\sf A}^{(1)}$ holds by construction.
 \medskip
 
 \noindent The following identities are used to show that
 \begin{align*}
 {\sf A}^{(1)2} {\sf A}^{*(2)} - (q^2+q^{-2}) {\sf A}^{(1)}{\sf A}^{*(2)} {\sf A}^{(1)} + {\sf A}^{*(2)} {\sf A}^{(1)2} = - H^2 (q^2-q^{-2})^2 {\sf A}^{*(2)}
 \end{align*}
 holds at  $ \tilde {\sf B}(r,s,t,u)$:
 the identity \eqref{eq:rec1} and also
 \begin{align*}
0&=C(r,t,s,u)C(u,r-1,t,s+1)(2 \theta^*_{s+u} - (q^2+q^{-2})\theta^*_{s+u+1})\\
& \quad +
 C(u,r,t,s)C(r,t+1,s,u-1)(2 \theta^*_{s+u} - (q^2+q^{-2}) \theta^*_{s+u-1}),\\
 0&=C(s,u,r,t)C(t,s-1,u,r+1)(2\theta^*_{s+u} - (q^2+q^{-2}) \theta^*_{s+u-1}) \\
 &\quad+ C(t,s,u,r)C(s,u+1,r,t-1)(2 \theta^*_{s+u} - (q^2+q^{-2}) \theta^*_{s+u+1}),\\
 - H^2 (q^2-q^{-2})^2 \theta^*_{s+u} 
 &= C(r,t,s,u)C(s+1,u,r-1,t)(2 \theta^*_{s+u} - (q^2+q^{-2}) \theta^*_{s+u+1}) \\
 &\quad+ C(t,s,u,r)C(u+1,r,t-1,s) (2 \theta^*_{s+u} - (q^2+q^{-2}) \theta^*_{s+u+1}) \\
 &\quad+ C(s,u,r,t)C(r+1,t,s-1,u)(2 \theta^*_{s+u} - (q^2+q^{-2}) \theta^*_{s+u-1}) \\
 &\quad + C(u,r,t,s)C(t+1,s,u-1,r)(2 \theta^*_{s+u} - (q^2+q^{-2}) \theta^*_{s+u-1}).
 \end{align*}
The identity \eqref{eq:rec2} is used to show that
 \begin{align*}
 {\sf A}^{*(2)2} {\sf A}^{(1)} - (q^2+q^{-2}) {\sf A}^{*(2)}{\sf A}^{(1)} {\sf A}^{*(2)} + {\sf A}^{(1)} {\sf A}^{*(2)2} = - H^2 (q^2-q^{-2})^2 {\sf A}^{(1)}.
 \end{align*}
\noindent The following identities are used to show that
\begin{align*}
\lbrack {\sf A}^{(1)}, \lbrack {\sf A}^{*(3)}, {\sf A}^{(2)}\rbrack_q \rbrack_q = \lbrack {\sf A}^{*(1)}, \lbrack {\sf A}^{(3)}, {\sf A}^{*(2)} \rbrack_q \rbrack_q
\end{align*}
holds at  $ \tilde {\sf B}(r,s,t,u)$:
\begin{align*}
0 &=C(r,u,t,s)C(r-1,t+1,s,u)(q^2\theta^*_{s+t+1}-\theta^*_{s+t})\\
 & \quad +C(r,t,s,u)C(r-1,u,t,s+1)(q^{-2}\theta^*_{s+t+1}-\theta^*_{s+t+2}), \\
 0 &= C(t,s,r,u)C(s,u,r+1,t-1)(q^2 \theta^*_{s+t-1}-\theta^*_{s+t}) \\
       &\quad + C(s,u,r,t)C(t,s-1,r+1,u)(q^{-2} \theta^*_{s+t-1}-\theta^*_{s+t-2}), \\
  0 &= C(s,r,u,t)C(s-1,u+1,r,t)(q^2 \theta^*_{s+t-1}-\theta^*_{s+t}) \\
  & \quad + C(s,u,r,t) C(s-1,r+1,u,t)(q^{-2} \theta^*_{s+t-1}-\theta^*_{s+t-2}), \\
  0 &= C(u,t,s,r)C(r,t,s+1,u-1)(q^2 \theta^*_{s+t+1}-\theta^*_{s+t}) \\
  &\quad + C(r,t,s,u)C(u,t,s+1,r-1)(q^{-2} \theta^*_{s+t+1}-\theta^*_{s+t+2})
   \end{align*}
  and also
\begin{align*}
  0 &= C(t,s,r,u)C(t-1,s,u,r+1)(q^2 \theta^*_{s+t-1}-\theta^*_{s+t}) \\
  &\quad + C(t,s,u,r)C(t-1,s,r,u+1)(q^{-2} \theta^*_{s+t-1}-\theta^*_{s+t-2}), \\
  0 &= C(r,u,t,s)C(u,r-1,t+1,s)(q^2 \theta^*_{s+t+1}-\theta^*_{s+t}) \\
  &\quad + C(u,r,t,s)C(r,u-1,t+1,s)(q^{-2} \theta^*_{s+t+1}-\theta^*_{s+t+2}), \\
  0 &= C(u,t,s,r)C(u-1,r,t,s+1)(q^2 \theta^*_{s+t+1}-\theta^*_{s+t}) \\
  &\quad + C(u,r,t,s)C(u-1,t+1,s,r)(q^{-2} \theta^*_{s+t+1}-\theta^*_{s+t+2}),\\
   0 &= C(s,r,u,t)C(t,s-1,u+1,r)(q^2 \theta^*_{s+t-1}-\theta^*_{s+t}) \\
   &\quad + C(t,s,u,r)C(s,r,u+1,t-1)(q^{-2} \theta^*_{s+t-1}-\theta^*_{s+t-2})
\end{align*}
and also
\begin{align*}
& C(r,s,u,t)(q \theta^*_{t+u+1}-q^{-1} \theta^*_{t+u})(q\theta^*_{u+s}-q^{-1} \theta^*_{u+s+1})\\
&= C(r,u,t,s)C(t+1,s,u,r-1)(q^2 \theta^*_{s+t+1}-\theta^*_{s+t}) \\
&\quad+ C(s,r,u,t)C(r,t,s-1,u+1)(q^2\theta^*_{s+t-1}-\theta^*_{s+t}) \\
&\quad+ C(r,t,s,u)C(s+1,r-1,u,t)(q^{-2} \theta^*_{s+t+1}-\theta^*_{s+t}) \\
&\quad + C(t,s,u,r)C(r,u+1,t-1,s)(q^{-2} \theta^*_{s+t-1}-\theta^*_{s+t}), \\
& C(t,r,s,u)(q \theta^*_{t+u-1}-q^{-1} \theta^*_{t+u})(q\theta^*_{u+s}-q^{-1} \theta^*_{u+s+1})\\
&= C(t,s,r,u)C(r+1,t-1,s,u)(q^2 \theta^*_{s+t-1}-\theta^*_{s+t}) \\
&\quad+ C(u,t,s,r)C(t,s+1,u-1,r)(q^2\theta^*_{s+t+1}-\theta^*_{s+t}) \\
&\quad+ C(r,t,s,u)C(t,s+1,r-1,u)(q^{-2} \theta^*_{s+t+1}-\theta^*_{s+t}) \\
&\quad + C(t,s,u,r)C(u+1,t-1,s,r)(q^{-2} \theta^*_{s+t-1}-\theta^*_{s+t})
\end{align*}
 and also
 \begin{align*}
 & C(s,u,t,r)(q \theta^*_{t+u+1}-q^{-1} \theta^*_{t+u})(q\theta^*_{u+s}-q^{-1} \theta^*_{u+s-1})\\
&= C(r,u,t,s)C(s,u,r-1,t+1)(q^2 \theta^*_{s+t+1}-\theta^*_{s+t}) \\
&\quad+ C(s,r,u,t)C(u+1,r,t,s-1)(q^2\theta^*_{s+t-1}-\theta^*_{s+t}) \\
&\quad+ C(s,u,r,t)C(r+1,u,t,s-1)(q^{-2} \theta^*_{s+t-1}-\theta^*_{s+t}) \\
&\quad + C(u,r,t,s)C(s,r,u-1,t+1)(q^{-2} \theta^*_{s+t+1}-\theta^*_{s+t}), \\
 & C(u,t,r,s)(q \theta^*_{t+u-1}-q^{-1} \theta^*_{t+u})(q\theta^*_{u+s}-q^{-1} \theta^*_{u+s-1})\\
&= C(t,s,r,u)C(u,r+1,t-1,s)(q^2 \theta^*_{s+t-1}-\theta^*_{s+t}) \\
&\quad+ C(u,t,s,r)C(s+1,u-1,r,t) (q^2\theta^*_{s+t+1}-\theta^*_{s+t}) \\
&\quad+ C(s,u,r,t)C(u,t,s-1,r+1)(q^{-2} \theta^*_{s+t-1}-\theta^*_{s+t}) \\
&\quad + C(u,r,t,s)C(t+1,s,r,u-1)(q^{-2} \theta^*_{s+t+1}-\theta^*_{s+t}).
\end{align*}

 \noindent We have verified some of the relations in Proposition \ref{lem:Nrels}. 
 The remaining relations in Proposition \ref{lem:Nrels} are verified  using the $S_3$-symmetry.
 


\section{Acknowledgements} The author thanks Kazumasa Nomura for reading the manuscript carefully and offering
valuable suggestions.

\bigskip


\noindent Paul Terwilliger \hfil\break
\noindent Department of Mathematics \hfil\break
\noindent University of Wisconsin \hfil\break
\noindent 480 Lincoln Drive \hfil\break
\noindent Madison, WI 53706-1388 USA \hfil\break
\noindent email: {\tt terwilli@math.wisc.edu }\hfil\break

\section{Statements and Declarations}

\noindent {\bf Funding}: The author declares that no funds, grants, or other support were received during the preparation of this manuscript.
\medskip

\noindent  {\bf Competing interests}:  The author  has no relevant financial or non-financial interests to disclose.
\medskip

\noindent {\bf Data availability}: All data generated or analyzed during this study are included in this published article.

\end{document}